\documentclass{amsart}

\usepackage{amssymb}
\usepackage{eepic}
\usepackage{color}

\allowdisplaybreaks

\numberwithin{equation}{section}

\newtheorem{theorem}{Theorem}[section]
\newtheorem{lemma}[theorem]{Lemma}
\newtheorem{proposition}[theorem]{Proposition}
\newtheorem{corollary}[theorem]{Corollary}
\newtheorem{remark}[theorem]{Remark}

\newtheorem{theoA}{Theorem}
\newtheorem{theoB}{Theorem}
\newtheorem{corA}{Corollary}
\newtheorem{corB}{Corollary}

\newtheorem{Rosenthal}{Rosenthal's theorem \cite{Ro}}
\newtheorem{Maurey}{Maurey's factorization theorem \cite{Mau}}

\theoremstyle{definition}

\newcommand{\N}{\mathbb{N}}
\newcommand{\Z}{\mathbb{Z}}
\newcommand{\R}{\mathbb{R}}
\newcommand{\C}{\mathbb{C}}

\newcommand{\Rad}{\operatorname{Rad}}

\newcommand{\E}{\mathcal{E}}

\newcommand{\U}{\mathcal{U}}

\newcommand{\la}{\mathbb{\lambda}}
\newcommand{\al}{\mathbb{\alpha}}

\newcommand{\ten}{\otimes}
\newcommand{\Om}{\Omega}
\newcommand{\eps}{\varepsilon}

\newcommand{\dem}{\noindent {\bf Proof. }}
\newcommand{\demA}{\noindent {\bf Proof of Theorem A. }}
\newcommand{\demB}{\noindent {\bf Proof of Theorem B. }}
\newcommand{\demCA}{\noindent {\bf Proof of Corollary A1. }}
\newcommand{\demCB}{\noindent {\bf Proof of Corollary A2. }}
\newcommand{\fin}{\hspace*{\fill} $\square$ \vskip0.2cm}
\newcommand{\prodd}{\prod\nolimits}
\newcommand{\limm}{\lim\nolimits}

\newcommand{\summ}{\sum\nolimits}

\begin{document}

\title[Noncommutative Maurey's factorization] {Maurey's
factorization theory \\ for operator spaces}

\author[Junge and Parcet]
{Marius Junge and Javier Parcet}


\footnote{Junge is partially supported by the NSF DMS-0556120.}
\footnote{Parcet is partially supported by \lq Programa Ram\'on y
Cajal 2005\rq${}$ \\ \null \hskip11.5pt and by Grants
MTM2007-60952 and CCG07-UAM/ESP-1664, Spain.}

\maketitle

\tableofcontents

\addtolength{\parskip}{+1ex}

\section*{Introduction}

In Banach space theory probabilistic techniques play a central
role. For example in the local theory of Banach spaces, geometric
properties of finite dimensional subspaces are proved from
probabilistic inequalities. The probabilistic approach not only
enriched Banach space theory, but also introduced Banach space
techniques in other areas such as probability or convex geometry.
A famous instance of such interplay is Maurey/Pisier's theory of
type and cotype. Their results are certainly inspired by
Rosenthal's work on subspaces of $L_p$. On the other hand, the
latter is strongly influenced by Grothendieck's notion of
absolutely summing maps, extended by Pietsch to $p>1$ and further
developed by Lindenstrauss/Pelczynski in their fundamental work on
Grothendieck's inequality.


All attempts to develop a similar theory for  operator spaces have
had only a limited success, so far. This is probably due to the
fact that there are many, if not too many, different operator
space structures on any Hilbert space. Indeed, in the local theory
of Banach spaces classification results typically measure the
distance of finite dimensional subspaces to Hilbert spaces and
then study critical indices, such as the best type $p$ or cotype
$q$ index \cite{Kr,MaP}. Therefore, the best one can hope for is
that for a given operator space there is a Hilbertian structure
which allows a similar local theory in the context of operator
spaces. A good illustration of this approach is Pisier's version
of Dvoretzky's theorem for operator spaces \cite{P3}. We will take
a different approach here.

This paper is inspired by the work on the \lq Grothendieck's
program\rq${}$ for operator spaces \cite{ER0,HM,J2,PS,X4}. To be
more precise, let us start by describing Rosenthal's theorem for
subspaces of $L_p$ and Maurey's factorization theorem. We first
recall some classical notions for a linear map $T: X \to Y$
between Banach spaces.
\begin{itemize}
\item[$\bullet$] $T$ has cotype $q$ if $$\Big( \sum_{k=1}^n
\|Tx_k\|_Y^q \Big)^{\frac1q} \ \le \ c_q(T) \, \Big( \mathbb{E}
\Big\| \sum_{k=1}^n \eps_k x_k \Big\|_X^q \Big)^{\frac1q},$$

\item[$\bullet$] $T$ is absolutely $(q,1)$-summing if $$\Big(
\sum_{k=1}^n \|Tx_k\|_Y^q \Big)^{\frac1q} \ \le \ \pi_{q,1}(T) \,
\sup_{\varepsilon_k = \pm 1} \Big\| \sum_{k=1}^n \eps_k x_k
\Big\|_X,$$

\item[$\bullet$] $T$ is $q$-summing if $$\Big( \sum_{k=1}^n
\|Tx_k\|_Y^q \Big)^{\frac1q} \ \le \ \pi_{q}(T) \,
\sup_{\|\phi\|_{X^*} \le 1} \Big( \sum_{k=1}^n | \langle \phi, x_k
\rangle |^q \Big)^{\frac1q}.$$
\end{itemize}
The constants $c_q(T), \pi_{q,1}(T), \pi_q(T)$ are the best ones
for which the inequalities hold.

\begin{Rosenthal}
Let $X \subset L_1$ be infinite dimensional and let $j: X \to L_1$
denote the inclusion map with adjoint $j^*: L_\infty \to X^*$.
Then, the following are equivalent\,$:$
\begin{itemize}
 \item[i)] $X$ embeds in $L_p$ for some $p>1$,
 \item[ii)] $X^*$ has cotype $q$ for some finite $q$,
 \item[iii)] $j^*$ is $(q,1)$-summing for some finite $q$.
\end{itemize}
\end{Rosenthal}

Using an adapted notion of $(q,1)$-concave maps, Rosenthal's
theorem remains true for infinite-dimensional subspaces of $L_p$
and $1 < p < 2$. The shortest way to prove Rosenthal's result is a
combination of the Grothendieck/Pietsch and Maurey's factorization
results. Indeed, Maurey's theorem (stated below) yields the hard
inclusion iii) $\Rightarrow$ i) in Rosenthal's result. The other
inclusions follows from well established facts in the theory.

\begin{Maurey}
Let $1 \le p < q < \infty$ and let $C(K)$ denote the space of
continuous functions in a compact Hausdorff space. Assume that the
linear map $T: C(K) \to X$ is absolutely $(p,1)$-summing. Then,
$T$ is $q$-summing and the following inequality holds
$$\pi_q(T) \le c(p,q) \, \pi_{p,1}(T).$$ This means that for any
absolutely $(p,1)$-summing map $T: C(K) \to X$, we may find a
probability measure $\mu$ and a linear map $w: L_q(K,\mu) \to X$
such that, if $j: C(K) \to L_q(K,\mu)$ denotes the natural
inclusion map, $T$ factorizes as $$T(x) = w \circ j(x).$$
\end{Maurey}

The main result of this paper is an operator space analog of
Maurey's theorem stated above and its natural generalization for
mappings $T: L_s \to X$. We refer to \cite{ER,P4} for basic
definitions on operator spaces. Motivated by Pisier's notion of a
completely $q$-summing operator \cite{P2}, we define a map $$T: X
\to Y$$ between operator spaces to be completely $(q,1)$-summing
if $$\pi_{q,1}^{cb}(T) = \big\| id \ten T: \ell_1 \ten_{\min} X
\to \ell_q(Y) \big\|_{cb} < \infty.$$ An expert in operator space
theory might think that it is more natural to take Schatten
classes $S_1$ and $S_q$ instead, see Remark \ref{Splpsumming}
below for a little discussion on this topic. Nevertheless, this
weaker notion is enough to obtain the operator space analog of
Maurey's factorization result.

\begin{theoA} \label{tha}
If $1 \le p < q < s < \infty$ and $X$ is an operator space, we
have\,$:$
\begin{itemize}
\item[i)] Let $A$ be a $C^*$-algebra and assume that the map $T: A
\to X$ is completely $(p,1)$-summing. Then, there exist positive
elements $\delta_1, \delta_2 \in L_{2q}(A^{**})$ and a map $w:
L_q(A^{**}) \to X$ such that $T(x) = w(\delta_1 x \delta_2)$ and
$$\|\delta_1\|_{2q} \|w\|_{cb} \|\delta_2\|_{2q} \le c(p,q) \,
\pi_{p,1}^{cb}(T).$$

\item[ii)] Let $\mathcal{M}$ be a von Neumann algebra and assume
that the map $T: \mathcal{M} \to X^*$ is normal and completely
$(p,1)$-summing. Then, there exist positive elements $d_1, d_2 \in
L_{2q}(\mathcal{M})$ and a map $v: L_q(\mathcal{M}) \to X^*$ such
that $T(x) = v(d_1 x d_2)$ and
$$\|d_1\|_{2q} \|v\|_{cb} \|d_2\|_{2q} \le c(p,q) \,
\pi_{p,1}^{cb}(T).$$

\item[iii)] Let $\mathcal{M}$ be a von Neumann algebra and assume
that the map $T: L_s(\mathcal{M}) \to X$ is completely
$(p,1)$-summing. Then, if $1/q = 1/s + 1/w$, there exist positive
elements $d_1, d_2 \in L_{2w}(\mathcal{M})$ and a completely
bounded map $v: L_q(\mathcal{M}) \to X$ such that $T(x) =
v(d_1xd_2)$ and $$\|d_1\|_{2w} \|v\|_{cb} \|d_2\|_{2w} \le
c(p,q,s) \, \pi_{p,1}^{cb}(T).$$
\end{itemize}
\end{theoA}

Note here that the analogue of a measure on $K$ is given by a
state $\phi$ on $A$. The natural analogue of the inclusion map
$id:C(K)\to L_p(K,\mu)$ is the positive map
$j_p(x)=d^{1/2p}xd^{1/2p}$ where $d$ is the positive density of
the state $\phi(x)=tr(dx)$ in $L_1(A^{**})$. Despite the analogy
of the results, a Banach space reader will have a hard time
recognizing similarities in the proof. The main difference relies
on the probabilistic part of the argument. Indeed, the new aspect
of the key embedding is based on our previous work \cite{JP3,JP4}.
Let us state it here since it might be of independent interest.
Let $X$ be an operator space and $\mathcal{M}$ be a von Neumann
algebra. Let us say that a linear map $T: X \to L_p(\mathcal{M})$
is $(p_1,p_2)$-convex if
$$k_{(p_1,p_2)}(T) = \big\| id \ten T: \ell_{p_1}(X) \to
L_p(\mathcal{M};\ell_{p_2}) \big\|_{cb} < \infty.$$

\begin{theoB}
Assume that $$T: X \to L_p(\mathcal{M})$$ is $(p_1,p_2)$-convex
and $1 \le p < q < (p_1 \wedge p_2) \le \infty$. Then, we have
$$\big\| T \ten id: S_q(X) \to L_p(\mathcal{M}; S_q) \big\|_{cb}
\le c(p,q,p_1,p_2) \, k_{(p_1,p_2)}(T).$$ 
\end{theoB}

We must emphasize that Theorems A and B hold for general von
Neumann algebras. The lack of a general theory of vector-valued
noncommutative $L_p$ spaces for arbitrary algebras forces us to
start with a careful analysis of the spaces we will handle along
the paper. Let us also note that in the special case $p=1$,
Theorem B is a dual version of Theorem A, and the corresponding
notion of concavity is even slightly weaker than the assumption in
Theorem A i). Our first application is of course an operator space
analog of Rosenthal's theorem. Our notion of cotype here will be
the following. Let $2 \le q < \infty$ and $$\Rad_q(X) = \Big\{
\summ_k \eps_k x_k \, \big| \ x_k \in X \Big\} \subset
L_q(\Sigma;X),$$ where the $\varepsilon_k$'s are independent $\pm
1$ Bernoulli's on a probability space $(\Sigma, \nu)$. Let $\iota$
be determined by $\iota(\varepsilon_k) = \delta_k$, where the
$\delta_k$'s form the canonical basis of $\ell_q$. Then we say
that a linear map $T: X \to Y$ between operator spaces has
cb-cotype $q$ if
$$c_q^{cb}(T) = \big\| \iota \ten T: \Rad_q(X) \to \ell_q(Y)
\big\|_{cb} < \infty.$$ An operator space $X$ has cb-cotype $q$ if
$id_X$ does. We refer to \cite{GP2,Lee,Pa0,Pa} for previous
attempts of defining a satisfactory notion of type and cotype for
operator spaces. In the following results, $p'$ will denote the
conjugate index of $p$, $\frac1p + \frac{1}{p'} = 1$. Rosenthal's
result takes the following form in the operator space setting.

\begin{corA}
If $1 \le p < 2$ and $X \subset L_p(\mathcal{M})$, t.f.a.e.
\begin{itemize}
\item[i)] There exists $p<q<2$ such that $X^*$ is of cb-cotype
$q'$.

\item[ii)] There exists $p<q<2$ such that $X^*$ is completely
$(q',1)$-summing.

\item[iii)] There exists $p<q<2$ such that $X$ completely embeds
into $L_q(\mathcal{M})$.
\end{itemize}
\end{corA}

In the category of Banach spaces, Rosenthal's theorem was recently
extended in \cite{JP2} for subspaces of noncommutative $L_p$
spaces. Although the relation with that result is obvious, we note
that Corollary A1 is not comparable since both hypotheses and
conclusions are stronger. Let us continue with the example of
Pisier's operator space $OH = [R,C]_{1/2}$. It is not too
difficult to prove that the identity map on $OH$ is completely
$(2,1)$ summing, see Lemma \ref{prel}. However, we know from
\cite{J2} that the strong version of the little Grothendieck
inequality fails
$$\mathcal{CB}(\mathcal{B(H)},OH) \not\subset
\Pi_2^o(\mathcal{B(H)},OH).$$

\begin{corB}
If $A$ is a $C^*$-algebra, $u: A \to OH$ is completely bounded if
and only if there exist positive elements $a,b \in L_1(A^{**})$
and a cb-map $w: L_p(A^{**}) \to OH$ for some $($\hskip-1pt
all\hskip1pt$)$ $2 < p < \infty$ such that $u = w (a^{1/2p} x
b^{1/2p})$. In particular, the isomorphism $\Pi_{p'}^o(OH,Y) =
\Pi_1^o(OH,Y)$ holds for $2 < p < \infty$ and any operator space
$Y$.
\end{corB}

We refer to \cite{P2} for the definition of the completely
$p$-summing norm $\pi_p^o$ and the space $\Pi_p^o(X,Y)$ of
completely $p$-summing maps $T: X \to Y$. This corollary vastly
improves on the results in \cite{JHL}. We see that $p>2$ is sharp
in this result, in contrast to what happens for Banach spaces. We
end up the paper with some further applications for Fourier
multipliers on discrete groups and other mappings between
noncommutative $L_p$ spaces.

\section{Vector-valued $L_p$ spaces} \label{S1}

Vector-valued, noncommutative $L_p$ spaces where introduced by
Pisier \cite{P2}. One of the main applications is a successful
understanding of noncommutative square and maximal functions. We
now discuss several settings for which vector-valued
noncommutative $L_p$ spaces are defined and which will be needed
below.

\subsection{The hyperfinite case}

In Pisier's setting, we assume that $\mathcal{M}$ is a hyperfinite
von Neumann algebra and $X$ is an arbitrary operator space. For
$1\le p<\infty$, the space $L_p(\mathcal{M}) = \lim_{\la}
L_p(\mathcal{M}_{\la})$ is a norm limit of finite dimensional von
Neumann subalgebras $\mathcal{M}_\lambda$. Therefore, it really
suffices to understand vector-valued Schatten $p$-classes. If
$R_p^m$ and $C_p^m$ stand for the row and column subspaces of
$S_p^m$, then define $$S_p^m(X) = C_p^m\ten_h X\ten_h  R_p^m.$$ In
operator space theory, the pairing $(a,b) = \mathrm{tr}(a^tb)$ is
chosen between $S_1^m$ and $M_m$. With respect to the paring
$\langle a,b \rangle = \mathrm{tr}(ab)$, we can reformulate the
main properties as follows:

\begin{itemize}
\item[a)] If $1 \le p \le\infty$, then
\begin{eqnarray*}
\mathrm{a1}) \quad \|x\|_{L_p(\mathcal{M};X)} & = & \inf_{x=ayb}
\|a\|_{2p} \|y\|_{\mathcal{M} \ten_{\min} X} \|b\|_{2p}, \\
\mathrm{a2}) \quad \|x\|_{L_p(\mathcal{M};X)} & = &
\sup_{\|a\|_{2p'},\|b\|_{2p'}\le 1} \,
\|axb\|_{L_1(\mathcal{M};X)}.
\end{eqnarray*}

\item[b)] If $1 \le p < \infty$, $L_p(\mathcal{M};X)^* =
L_{p'}(\mathcal{M}^{\mathrm{op}};X^*)$ with respect to the bracket
$$\Big\langle \summ_j a_j\ten x_j^*, \summ_k b_k \ten x_k
\Big\rangle = \summ_{j,k} \mathrm{tr}(a_j b_k) \, \langle x_j^*,
x_k \rangle.$$

\end{itemize}

\subsection{Amalgamated and conditional $L_p$ spaces}

Let us now recall some new noncommutative function spaces from
\cite{JP4} which will be essential below. Let $\mathcal{M}$ be an
arbitrary von Neumann algebra and let $\mathcal{R}$ stand for the
matrix amplification $\mathcal{M} \bar\ten \mathcal{B}(\ell_2)$.
In what follows we shall work with indices represented in the
following solid of $\R^3$ $$\mathsf{K} = \Big\{(1/u,1/v,1/q) \,
\big| \ 2 \le u,v \le \infty, \ 1 \le q \le \infty, \ 1/u + 1/q +
1/v \le 1 \Big\}.$$ Given $1 \le p \le \infty$ such that $\frac1p
= \frac1u + \frac1q + \frac1v$ for some $(\frac1u, \frac1v,
\frac1q)$ in $\mathsf{K}$, we define the corresponding amalgamated
$L_p$ space as the subspace
$L_{\underline{u}q\underline{v}}(\mathcal{R}; \mathcal{M})$ of
$L_p(\mathcal{R})$ equipped with the norm
$$\|x\|_{\underline{u}q\underline{v}} = \inf \Big\{
\|a\|_{L_u(\mathcal{M})} \|y\|_{L_q(\mathcal{R})}
\|b\|_{L_v(\mathcal{M})} \, \big| \ x = a y b \Big\}.$$ We shall
also be interested in the duals of amalgamated $L_p$ spaces. To
that end given $(\frac1u, \frac1v, \frac1p) \in \mathsf{K}$ and
$\frac1s = \frac1u + \frac1p + \frac1v$, we define the
corresponding conditional $L_p$ space
$L_{\overline{u}s\overline{v}}(\mathcal{R}; \mathcal{M})$ as the
completion of $L_p(\mathcal{R})$ with respect to the norm
$$\|x\|_{\overline{u}s\overline{v}} = \sup \Big\{ \|\alpha x
\beta\|_{L_s(\mathcal{R})} \, \big| \
\|\alpha\|_{L_u(\mathcal{M})}, \|\beta\|_{L_v(\mathcal{M})} \le 1
\Big\}.$$ We refer to \cite{JP4} for a more detailed exposition
and note in passing that we have changed/improved our terminology
for amalgamated and conditional $L_p$'s. Now we collect the main
complex interpolation and duality properties from \cite{JP4}. Let
$\mathsf{K}_0$ denote the interior of $\mathsf{K}$. Then we have:

\vskip3pt

\begin{itemize}
\item[i)] $L_{\underline{u}q\underline{v}}(\mathcal{R};
\mathcal{M})$ is a Banach space.

\vskip3pt

\item[ii)] $L_{\underline{u_{\theta}} q_{\theta}
\underline{v_{\theta}}}(\mathcal{R}; \mathcal{M})$ is
isometrically isomorphic to
$$\big[L_{\underline{u_0}q_0\underline{v_0}}(\mathcal{R};
\mathcal{M}), L_{\underline{u_1}q_1\underline{v_1}}(\mathcal{R};
\mathcal{M}) \big]_{\theta}^{\null},$$ with $(\frac{1}{u_\theta},
\frac{1}{q_\theta}, \frac{1}{v_\theta}) = (\frac{1-\theta}{u_0} +
\frac{\theta}{u_1}, \frac{1-\theta}{q_0} + \frac{\theta}{q_1},
\frac{1-\theta}{v_0} + \frac{\theta}{v_1})$.

\vskip3pt

\item[iii)] If $(1/u,1/v,1/q) \in \mathsf{K}_0$ and $1 - 1/p = 1/u
+ 1/q + 1/v$ $$\hskip20pt \big(
L_{\underline{u}q\underline{v}}(\mathcal{R}; \mathcal{M}) \big)^*
= L_{\overline{u}q'\overline{v}}(\mathcal{R}; \mathcal{M}) \quad
\mbox{and} \quad \big( L_{\overline{u}q'\overline{v}}(\mathcal{R};
\mathcal{M}) \big)^* =
L_{\underline{u}q\underline{v}}(\mathcal{R}; \mathcal{M}),$$ with
respect to the antilinear duality bracket $\langle x,y \rangle =
\mathrm{tr}(x^*y)$.
\end{itemize}

\subsection{Mixed norms I}

\label{MnI}

The definition of amalgamated and conditional $L_p$ spaces was
mainly inspired by Pisier's fundamental identities for the mixed
norm spaces $L_p(\mathcal{M}_1; L_q(\mathcal{M}_2))$ with
$\mathcal{M}_1$ hyperfinite. Given $1 \le p,q \le \infty$ and
$\frac1r = |\frac1p - \frac1q|$, we have $$\|x\|_{L_p(L_q)} =
\left\{ \begin{array}{ll} \inf \Big\{
\|\alpha\|_{L_{2r}(\mathcal{M}_1)} \|y\|_{L_q(\mathcal{M}_1
\bar\otimes \mathcal{M}_2)} \|\beta\|_{L_{2r}(\mathcal{M}_1)} \,
\big| \ x = \alpha y \beta \Big\} & \mathrm{if} \ p \le q, \\
[8pt] \sup \Big\{ \|\alpha x \beta \|_{L_q(\mathcal{M}_1
\bar\otimes \mathcal{M}_2)} \, \big| \,
\|\alpha\|_{L_{2r}(\mathcal{M}_1)},
\|\beta\|_{L_{2r}(\mathcal{M}_1)} \le 1 \Big\} & \mathrm{if} \ p
\ge q. \end{array} \right.$$ The extension to arbitrary von
Neumann algebras is a matter of regarding these spaces as
amalgamated and conditional $L_p$ spaces. Indeed, given any von
Neumann algebra $\mathcal{M}$ and $\mathcal{R} = \mathcal{M}
\bar\ten \mathcal{B}(\ell_2)$, we may define
$$L_p(\mathcal{M}; S_q) = \left\{
\begin{array}{ll} L_{\underline{2r}q\underline{2r}}(\mathcal{R};
\mathcal{M}) & \mathrm{if} \ p \le q, \\ [3pt]
L_{\overline{2r}q\overline{2r}}(\mathcal{R}; \mathcal{M}) &
\mathrm{if} \ p \ge q, \end{array} \right.$$

\begin{remark}
\emph{We define an operator space structure on $L_p(\mathcal{M};
S_q)$ by complex interpolation. It thus suffices to provide the
o.s.s. of the endpoints $L_p(\mathcal{M}; S_q)$ for $p,q \in
\{1,\infty\}$. If $p=q$ the definition is obvious, while
$L_1(\mathcal{M}; S_\infty)$ embeds into the dual of
$L_\infty(\mathcal{M}; S_1)$. Hence, it just remains to understand
the o.s.s. of the latter one. According to \cite{P4}, we may
define the operator space $L_\infty(\mathcal{M}; S_1)$ as the
quotient
$$\mathcal{M} \otimes_h S_1 \otimes_h \mathcal{M} \big/
\mathrm{ker} \, q,$$ by the quotient map $q(a \ten x \ten b) = ab
\otimes x$. Moreover, we also find}
\begin{itemize}
\item \emph{Complex interpolation also gives}
$S_p^n(L_p(\mathcal{M};S_q)) = L_p(M_n \otimes \mathcal{M}; S_q)$.

\item \emph{The same argument provides an o.s.s. for $A(\ell_1)$
with $A$ any $C^*$-algebra.}
\end{itemize}
\end{remark}

\noindent Then we easily find that
\begin{itemize}
\item[a)] If $1 \le p \le\infty$, then
\begin{eqnarray*}
\mathrm{a1}) \quad \|x\|_{L_p(\mathcal{M};S_q^m)} & = &
\inf_{x=ayb}
\|a\|_{2p} \|y\|_{L_\infty(\mathcal{M}; S_q^m)} \|b\|_{2p}, \\
\mathrm{a2}) \quad \|x\|_{L_p(\mathcal{M};S_q^m)} & = &
\sup_{\|a\|_{2p'},\|b\|_{2p'}\le 1} \,
\|axb\|_{L_1(\mathcal{M};S_q^m)}.
\end{eqnarray*}

\item[b)] If $1 \le p < \infty$, we have $L_p(\mathcal{M};S_q^m)^*
= L_{p'}(\mathcal{M};S_{q'}^m)$.
\end{itemize}
Note that $\mathcal{M} \ten_\mathrm{min} X$ in the hyperfinite
case is replaced here by $L_\infty(\mathcal{M};S_q^m)$. It should
be noticed that we still have $L_p(\mathcal{M};S_q^m) =
[L_p(\mathcal{M}; S_\infty^m), L_p(\mathcal{M},S_1^m)]_{1/q}$ for
general von Neumann algebras, see \cite{JX4,JXnew}. It is well
known since \cite{J1} that the norms of the boundary points are
given by
\begin{eqnarray*}
\|x\|_{L_p(\mathcal{M}; S_\infty^m)} & = & \inf_{x = a y b} \
\|a\|_{L_{2p}(\mathcal{M})} \|y\|_{M_m(\mathcal{M})}
\|b\|_{L_{2p}(\mathcal{M})}, \\ \|x\|_{L_p(\mathcal{M}; S_1^m)} &
= & \inf_{x_{ij} = \stackrel{}{\sum_k} a_{ik}b_{jk}} \Big\| \Big(
\summ_{i,k} a_{ik}a_{ik}^* \Big)^{\frac12} \Big\|_{2p} \, \Big\|
\Big( \summ_{j,k} b_{jk}^* b_{jk} \Big)^{\frac12} \Big\|_{2p}.
\end{eqnarray*}

\begin{remark}
\emph{We have just considered the amplification algebra
$\mathcal{R}$ since we shall be mainly interested in mixed-norms
with values in matrix algebras. Nevertheless at some points in
this paper we will handle spaces of the form $L_p(\mathcal{M}_1;
L_q(\mathcal{M}_2))$ for non-hyperfinite $\mathcal{M}_j$. In this
case the notions of amalgamated and conditional $L_p$ spaces
remains unchanged, see \cite{JP4} for further details.}
\end{remark}

\subsection{Asymmetric Schatten classes}

Given any von Neumann algebra $\mathcal{M}$, we write
$L_2^r(\mathcal{M})$ and $L_2^c(\mathcal{M})$ to denote the
row/column quantizations on $L_2(\mathcal{M})$ and consider the
operator spaces $$L_u^r(\mathcal{M}) =
[\mathcal{M},L_2^r(\mathcal{M})]_\frac{2}{u} \quad \mbox{and}
\quad L_v^c(\mathcal{M}) =
[\mathcal{M},L_2^c(\mathcal{M})]_\frac{2}{v}.$$ In fact, a
rigorous definition should take Kosaki's embeddings into account
as done in \cite[Identity (1.3)]{JP4}, but we shall ignore such
formalities here. We have the complete isometry $L_p(\mathcal{M})
= L_{2p}^r(\mathcal{M}) \ten_{\mathcal{M},h}
L_{2p}^c(\mathcal{M})$, where $\otimes_{\mathcal{M},h}$ stands for
the $\mathcal{M}$-amalgamated Haagerup tensor product. This
motivates the definition of the asymmetric spaces
$$L_{(2u,2v)}(\mathcal{M}) = L_{2u}^r(\mathcal{M})
\ten_{\mathcal{M},h} L_{2v}^c(\mathcal{M}) =
L_{\underline{2u}\infty\underline{2v}}(\mathcal{M};
\mathcal{M}).$$ These spaces were originally defined in \cite{JP1}
for finite matrix algebras, where the definition simplifies in
terms of ordinary Haagerup tensors. In this case, given an
arbitrary operator space $X$, we may as well consider the
vector-valued space as $$S_{(2u,2v)}^m(X) = C_u^m \ten_h X \ten_h
R_v^m.$$ Its module behavior is explained better by
$$S_{(2u,2v)}^m(X) = S_{(2u,\infty)}^m \otimes_{M_m,h}
S_\infty^m(X) \otimes_{M_m,h} S_{(\infty,2v)}^m.$$ Again by
interpolation, we find a natural o.s.s. for
$L_{(2u,2v)}(\mathcal{M})$ and we see that
$$C_u^n \ten_h L_{(2u,2v)}(\mathcal{M}) \ten_h R_v^n =
L_{(2u,2v)}(M_n \ten \mathcal{M}).$$ According to \cite{JP1}, we
have
\begin{itemize}
\item[a)] If $1 \le u,v \le\infty$, then
\begin{eqnarray*}
\mathrm{a1}) \quad \|x\|_{S_{(2u,2v)}^m(X)} & = & \inf_{x=ayb}
\|a\|_{2u} \|y\|_{S_\infty^m(X)} \|b\|_{2v}, \\
\mathrm{a2}) \quad \|x\|_{S_{(2u,2v)}^m(X)} & = &
\sup_{\|a\|_{2u'},\|b\|_{2v'}\le 1} \, \|axb\|_{S_1^m(X)}.
\end{eqnarray*}

\item[b)] If $1 \le u,v \le \infty$, $S_{(2u,2v)}^m(X)^* =
S_{(2u',2v')}^m(X^*)$ with respect to the bracket
$$\Big\langle \summ_j a_j \ten x_j^*, \summ_k b_k \ten x_k
\Big\rangle = \summ_{j,k} \mathrm{tr}(a_j b_k) \, \langle x_j^*,
x_k \rangle.$$

\end{itemize}

\subsection{Mixed norms II}

The next family of spaces are noncommutative $L_p$ spaces with
values in asymmetric Schatten classes. Namely, let us recall the
spaces $L_\infty(\mathcal{M}; C_q) = [L_\infty(\mathcal{M};
C_\infty), L_\infty(\mathcal{M}; R_\infty)]_{1/q}$ defined from
the row/column spaces $L_\infty(\mathcal{M}; R_\infty) =
\mathcal{M} \bar\ten R$ and $L_\infty(\mathcal{M}; C_\infty) = C
\bar\ten \mathcal{M}$. The spaces $L_\infty(\mathcal{M}; C_q)$
were already considered by Pisier for semifinite von Neumann
algebras \cite{P5} and by Haagerup for general von Neumann
algebras \cite{H4}. Define
\begin{eqnarray*}
L_{2p}^r(\mathcal{M}; C_q) & = & L_{2p}^r(\mathcal{M})
\ten_{\mathcal{M},h} L_\infty(\mathcal{M};C_q),
\\ L_{2p}^c(\mathcal{M}; R_q) & = & L_\infty(\mathcal{M};R_q)
\ten_{\mathcal{M},h} L_{2p}^c(\mathcal{M}).
\end{eqnarray*}
These spaces satisfy:
\begin{itemize}
\item[i)] $L_{2p}^r(\mathcal{M}; C_p) =
C_{2p}(L_{2p}(\mathcal{M}))$ isometrically.

\item[ii)] $L_{2p_\theta}^r(\mathcal{M};C_{q_\theta}) = \big[
L_{2p_0}^r(\mathcal{M};C_{q_0}), L_{2p_1}^r(\mathcal{M};C_{q_1})
\big]_\theta$ isometrically.
\end{itemize}
and analogous properties hold for the adjoint spaces. Indeed, the
second property follows from a nowadays standard interpolation
technique originated in \cite{P5} and further developed in
\cite{JP4,X}. The first property is clear for $p = \infty$ and it
then suffices by interpolation to consider the case $p=1$. Again,
this is standard by applying Pisier's factorization trick in
\cite{P2}. If $\frac1s = |\frac1p - \frac1q|$, the norm can
written as follows $$\|x\|_{L_{2p}^r(\mathcal{M}; C_q)} = \left\{
\begin{array}{ll} \inf_{x=\alpha y} \|\alpha\|_{2s}
\|y\|_{C_{2q}(L_{2q}(\mathcal{M}))} & \mbox{if} \ p \le q, \\
[5pt] \sup_{\|\alpha\|_{2s} \le 1} \hskip3pt \|\alpha
x\|_{C_{2q}(L_{2q}(\mathcal{M}))} & \mbox{if} \ p \ge q.
\end{array} \right.$$ If $1 \le p \le q_1 \wedge q_2 \le \infty$
and $1 \le s_1, s_2 \le \infty$ satisfy $\frac{1}{s_j} = \frac1p -
\frac{1}{q_j}$ and $\mathcal{R} = \mathcal{M} \bar\ten
\mathcal{B}(\ell_2)$, we set $$L_p(\mathcal{M}; C_{q_1} \ten_h
R_{q_2}) = L_{2p}^r(\mathcal{M}; C_{q_1}) \ten_{\mathcal{M},h}
L_{2p}^c(\mathcal{M}; R_{q_2}).$$ By interpolation, it is
compatible with the symmetric case given in Paragraph \ref{MnI}.

\begin{remark}
\emph{Connes' characterization of hyperfiniteness can be rephrased
by the condition $L_\infty(\mathcal{M}; OH) \simeq \mathcal{M}
\ten_{\min} OH$, see \cite{P1}. Thus in general we have to accept
that the norms considered so far are different, but consistent.
Namely, we have seen that the asymmetric/nonhyperfinite norms
generalize the symmetric/hyperfinite ones respectively. Thus,
there should be no ambiguity of what definition is being used
along the text.}
\end{remark}

\begin{remark} \label{linn}
{\rm If $A$ is a $C^*$-algebra, consider the norm
\[ \|(x_k)\|_{A(\ell_1)} = \inf_{x_k = \sum_j a_{kj}b_{kj}} \Big\|
\Big( \summ_{j,k} a_{jk} a_{kj}^* \Big)^{\frac12} \Big\| \Big\|
\Big( \summ_{j,k} b_{jk}^* b_{kj} \Big)^{\frac12} \Big\| .\]
Replacing $A$ by $M_m(A)$, we see that $A(\ell_1) \subset \ell_1
\ten_\mathrm{min} A$ is a complete contraction and according to an
unpublished work of Haagerup, this is an isometry only for
$C^*$-algebras with Lance's weak expectation property. At any
rate, we see that every completely $(p,1)$-summing map $T:A\to X$
satisfies $$\big\| id \ten T: A(\ell_1) \to \ell_p(X) \big\|_{cb}
\le \pi_{(p,1)}^{cb}(T).$$ Indeed, for Theorem A i) only this
weaker assumption of $(p,1)$-concavity is required. This concavity
is the cb-version of Pisier's notion of $(p,1)$ $C^*$-summability
in \cite{P0}.}
\end{remark}

\begin{remark}
\emph{In the hyperfinite and asymmetric cases we considered
arbitrary operator spaces and specific von Neumann algebras. In
the mixed-norm cases the situation is the opposite. There exists
an intermediate notion of $L_p(\mathcal{M}; X)$ valid for QWEP von
Neumann algebras and operator spaces contained in any
$C^*$-algebra with the local lifting property. The notion was
developed in \cite{J4} and is based on the hyperfinite theory
replacing norm approximations by ultraproducts. Some arguments in
this paper could be slightly simplified if we restricted to work
over QWEP von Neumann algebras.}
\end{remark}

\section{Key probabilistic estimates}
\label{S2}

In this section we use the theory of vector-valued $L_p$ spaces in
connection with convexifying operators. This leads to a change of
density which will be crucial for our proof of Maurey's theorem.

\subsection{A cb-embedding for $S_q(X)$} Let us consider a weight
function $w$ indexed over the integers $\Z$ and define the Hilbert
space $\ell_2(w)$ determined by the following norm
$$\Big\| \summ_n w_n \delta_n \Big\|_{\ell_2(w)} = \Big(
\summ_n w_n |\alpha_n|^2 \Big)^{\frac12}.$$ If $\ell_2^r(w)$ and
$\ell_2^c(w)$ denote the row/column o.s.s. on $\ell_2(w)$, we set
$$\ell_2^{r_p}(w) = \big[ \ell_2^r(w),
\ell_2^c(w) \big]_\frac1p \quad \mbox{and} \quad \ell_2^{c_p}(w) =
\big[ \ell_2^c(w), \ell_2^r(w) \big]_\frac1p.$$ Most of the time,
our weights will be of the form $w_n = \lambda^n$ for some
$\lambda > 1$. In that cases we will write $w_\lambda$ and
$\ell_2(w_\lambda)$ instead. Our first step will be a description
of $S_q(X)$ closely related to Xu's characterization \cite{X5} of
$R_q$ and $C_q$. Although it also follows from a more general
argument in \cite{JX5}, we give here a concrete approach for
completeness. In what follows we shall write $\alpha \lesssim
\beta$ to indicate the existence of an absolute constant $c$ such
that $\alpha \le c \beta$. We begin with a well-known observation
comparing the $J$ and $K$ methods as in \cite{J2,X5}.

\begin{lemma} \label{dis1}
Let $A$ and $B$ be non-singular positive operators on a Hilbert
space $\mathcal{H}$ and assume further than $A$ and $B$ commute.
If $0<\theta<1$ and $\la > 1$, let us consider the constants
$$c_1(\la,\theta) = \sqrt{\frac{1}{\la^{\theta}-1}
+ \frac{1}{\la^{1-\theta}-1}} \quad , \quad c_2(\la,\theta) =
\sqrt{\frac{1}{1 - \la^{-\theta}} + \frac{1^{\null}}{1 -
\la^{-(1-\theta)}}}.$$ Then, the equivalence $c_1(\lambda,\theta)
\, \alpha \lesssim \beta \lesssim c_2(\lambda,\theta) \, \alpha$
holds with
\begin{eqnarray*}
\alpha & = & \big\| A^\theta B^{1-\theta}x \big\|_\mathcal{H}, \\
\beta & = & \inf_{x = y_n + z_n} \Big( \sum_{n \in \Z}
\la^{n(1-\theta)} \|Ay_n\|_\mathcal{H}^2 \Big)^\frac12 + \Big(
\sum_{n \in \Z} \la^{-n\theta} \|Bz_n\|_\mathcal{H}^2
\Big)^\frac12.
\end{eqnarray*}
A duality argument also gives $c_1(\lambda,\theta) \, \beta'
\lesssim \alpha' \lesssim c_2(\lambda,\theta) \, \beta'$ with
\begin{eqnarray*}
\alpha' & = & \big\| A^{1-\theta} B^{\theta}x \big\|_\mathcal{H},
\\ \beta' & = & \inf_{A^{1-2\theta}B^{2\theta-1}x = \sum_n z_n} \Big(
\sum_{n \in \Z} \la^{-n(1-\theta)} \|Az_n\|_\mathcal{H}2
\Big)^\frac12 + \Big( \sum_{n \in \Z} \la^{n\theta}
\|Bz_n\|_\mathcal{H}^2 \Big)^\frac12.
\end{eqnarray*}
\end{lemma}

\dem By simultaneous diagonalization, it suffices to prove the
first assertion for diagonal operators $A = D_a$ and $B = D_b$.
Then it is clear that the term $\beta$ is equivalent to
$$\inf_{x_k = y_{nk} + z_{nk}}
\Big( \sum_{n \in \Z} \la^{n(1-\theta)} \sum_{k \in \N} |a_k|^2
|y_{nk}|^2 + \sum_{n \in \Z} \la^{-n\theta} \sum_{k \in \N}
|b_k|^2 |z_{nk}|^2 \Big)^\frac12.$$ This equals $$\Big( \sum_{k
\in \N} \, \inf_{x_k = y_{nk} + z_{nk}} \big[ \sum_{n \in \Z}
\la^{n(1-\theta)} |a_k y_{nk}|^2 + \sum_{n \in \Z} \la^{-n\theta}
|b_k z_{nk}|^2 \big] \Big)^\frac12.$$ Thus, it suffices to prove
the assertion for $k$ fixed and then we may even assume that
$x_k=1$ by normalization. This reduces the assertion to scalars
and we therefore claim that
$$c_1(\lambda,\theta) \, a^\theta b^{1-\theta} \lesssim
\mathcal{B} \lesssim c_2(\lambda,\theta) \, a^\theta
b^{1-\theta}$$ for $a,b > 0$ and
$$\mathcal{B} = \inf_{1=\gamma_n+\rho_n} \Big( \sum_{n \in \Z}
\la^{n(1-\theta)} |a \gamma_n|^2 + \sum_{n \in \Z} \la^{-n\theta}
|b \rho_n|^2 \Big)^\frac12.$$ Let us start with an easy
observation
$$\inf_{1=\gamma+\rho} \delta |\gamma|^2 + \sigma |\rho|^2 =
\inf_{0 \le t \le 1} \delta t^2+ \sigma (1-t)^2 = \frac{\delta
\sigma}{\delta+\sigma} \sim \min(\delta,\sigma)$$ holds for all
$\delta, \sigma > 0$. Going back to our claim, since
$|\gamma_n|+|\rho_n| \ge 1$, it therefore suffices to consider
$\gamma_n$ and $\rho_n$ positive in the right hand side above.
This leads to the following estimate
\begin{eqnarray*}
\mathcal{B}^2 & = & \inf_{1=\gamma_n+\rho_n} \sum_{n \in \Z} \big[
a^2 \la^{n(1-\theta)} |\gamma_n|^2 + b^2 \la^{-n\theta} |\rho_n|^2
\big] \sim \sum_{n \in \Z} \min(a^2 \la^{n(1-\theta)}, b^2
\la^{-n\theta}) \\ & = & \sum_{\la^{-n} \ge a^2/b^2} a^2
\la^{n(1-\theta)} + \sum_{\la^{-n} < a^2/b^2} b^2 \la^{-n\theta} \ =
\ a^2 \frac{\la^{n_0(1-\theta)}}{1-\la^{-(1-\theta)}} + b^2
\frac{\la^{-(n_0+1)\theta}}{1-\la^{-\theta}}.
\end{eqnarray*}
Here $n_0$ is chosen so that $\la^{-(n_0+1)} < a^2/b^2 \le
\la^{-n_0}$ and this gives $$c_1(\la,\theta)^2 \, a^{2\theta}
b^{2(1-\theta)} \lesssim \mathcal{B}^2 \lesssim c_2(\la,\theta)^2 \,
a^{2\theta} b^{2(1-\theta)}.$$ Hence, the first assertion follows.
To prove the second assertion, given a positive non-singular
operator $L$ acting on $\mathcal{H}$, we denote by $\mathcal{H}_L$
the Hilbert space equipped with the norm
$$\|x\|_{\mathcal{H}_L} = \|Lx\|_\mathcal{H}.$$ By the
first assertion, we know that
$\mathcal{H}_{A^{\theta}B^{1-\theta}}$ is isomorphic (up to the
constants $c_j(\la,\theta)$) to the subspace of constant sequences
in $\ell_2(\la^{1-\theta};\mathcal{H}_A) +
\ell_2(\la^{-\theta};\mathcal{H}_B)$. Since
$\mathcal{H}_{A^{\theta}B^{1-\theta}}$ is a Hilbert space, it is
isometric to its dual. In particular, recalling that
$$A^{1-\theta}B^\theta x = A^\theta B^{1-\theta} \big( A^{1-2\theta}
B^{2\theta-1}x \big),$$ we find that its norm in $\mathcal{H}$ is
equivalent to the norm of $A^{1-2\theta}B^{2\theta-1}x$ in the
quotient of $\ell_2(\la^{-(1-\theta)};\mathcal{H}_A) \cap
\ell_2(\la^{\theta};\mathcal{H}_B)$ by the subspace of mean zero
sequences. Writing this down we obtain the second assertion. The
proof is complete. \fin

To continue, we need to introduce Xu's terminology in \cite{X5}.
We will only define the column spaces, but we shall freely use
below the row analogs which are defined in the obvious way. Let
$$\ell_2(w;\ell_2) = \ell_2(\Z,w; \ell_2(\N)) \quad
\mbox{with norm} \quad \Big\| \sum_{n,k} x_{nk} \otimes
\delta_{nk} \Big\| = \Big( \sum_{n \in \Z} w_n \sum_{k \in \N}
|x_{nk}|^2 \Big)^\frac12.$$ Define $\ell_2(w_\lambda;\ell_2)$
similarly, let $\ell_2^{c_p}(w_\lambda;\ell_2) =
[\ell_2^c(w_\lambda;\ell_2),
\ell_2^r(w_\lambda;\ell_2)]_{\frac1p}$ and set
$$\mathcal{G}_{c_p,c_q}^K(w_\la,\theta) =
\ell_2^{c_p}(w_\la^{-\theta};\ell_2) +
\ell_2^{c_q}(w_\la^{1-\theta};\ell_2) \quad \mbox{with} \quad
w_\lambda^\eta = w_{\lambda^\eta}.$$ Let
$\mathcal{C}_{c_p,c_q}^K(w_\la,\theta)$ denote the subspace of
$\Z$-constant sequences. Using the bracket
$$\big\langle (a_{nk}), (b_{nk}) \big\rangle = \sum_{n \in \Z}
\sum_{k \in \N} \overline{a}_{nk} b_{nk},$$ the dual spaces are
\begin{eqnarray*}
\big( \mathcal{G}_{c_p,c_q}^K(w_\la,\theta) \big)^* & = &
\mathcal{G}_{c_{p'},c_{q'}}^J(w_\la^{-1},\theta), \\ \big(
\mathcal{C}_{c_p,c_q}^K(w_\la,\theta) \big)^* & = &
\mathcal{C}_{c_{p'},c_{q'}}^J(w_\la^{-1},\theta),
\end{eqnarray*}
with operator space structures given by
\begin{eqnarray*}
\Big\| \sum_{n \in \Z} \sum_{k=1}^\infty x_{nk} \ten e_{(n,k),1}
\Big\|_{M_m(\mathcal{G}_J)} & = & \max \big\{ n_p(x), n_q(x) \big\}, \\
\Big\| \sum_{n \in \Z} \sum_{k=1}^\infty x_{nk} \ten e_{(n,k),1} +
\mathcal{C}_K ^\perp \Big\|_{M_m(\mathcal{C}_J)} & = &
\inf_{\sum_n x_{nk} - z_{nk} = 0} \max \big\{ n_p(z), n_q(z)
\big\}
\end{eqnarray*}
where the norms $n_p(\xi)$ and $n_q(\xi)$ are given by
\begin{eqnarray*}
n_p(\xi) & = & \Big\| \sum_{n \in \Z} \la^{n\theta/2}
\sum_{k=1}^\infty \xi_{nk} \ten e_{(n,k),1} \Big\|_{M_m(C_{p'})}, \\
n_q(\xi) & = & \Big\| \sum_{n \in \Z} \la^{-n(1-\theta)/2}
\sum_{k=1}^\infty \xi_{nk} \ten e_{(n,k),1} \Big\|_{M_m(C_{q'})}.
\end{eqnarray*}

The following result is closely related to \cite[Section 2]{X3}.
However  we have to review the argument in order to understand the
generalization presented below.

\begin{lemma} \label{Xu}
If $p_0<q<p_1$ with $\frac{1}{q} = \frac{1-\theta}{p_0} +
\frac{\theta}{p_1}$ and $\la > 1$, then $$R_q \simeq_{cb}
\mathcal{C}_{r_{p_0},r_{p_1}}^K(w_\la,\theta) \quad \mbox{and}
\quad C_q \simeq_{cb}
\mathcal{C}_{c_{p_0},c_{p_1}}^K(w_\la,\theta).$$ The constant of
these complete isomorphisms only depend on $\la$ and $\theta$.
\end{lemma}

\dem Since both cb-isomorphisms are proved in the same way, we
only argue with column spaces. Let us first show that the
inclusion $\mathcal{C}_{c_{p_0},c_{p_1}}^K(w_\la,\theta) \subset
C_q$ is completely bounded. We recall the o.s.s. of $C_q$ from the
main result in \cite{X}
$$\Big\| \sum_{k=1}^\infty x_k \ten e_{k,1} \Big\|_{M_m(C_q)}
= \sup_{\|a\|_{S_{2q}^m}, \|b\|_{S_{2q'}^m} \le 1} \Big(
\sum_{k=1}^\infty \|ax_kb\|_2^2 \Big)^\frac12.$$ We may clearly
assume that $a$ and $b$ are positive and invertible. Let us denote
by $\mathcal{L}_a(x)=a x$ and $\mathcal{R}_{b}(x)=xb$ the
left/right actions. We define
$A=\mathcal{L}_{a^{q/p_1}}\mathcal{R}_{b^{q'/p_1'}}$ and
$B=\mathcal{L}_{a^{q/p_0}}\mathcal{R}_{b^{q'/p_0'}}$. Then we
apply Lemma \ref{dis1} to $x = \sum_k x_k \ten e_{k,1}$ and deduce
that we have $$\Big( \sum_{k=1}^\infty \|ax_kb\|_2^2 \Big)^\frac12
\ = \ \big\| A^{\theta} B^{1-\theta} x \big\|_2 \ \lesssim \
c_1(\la,\theta)^{-1} \inf_{x_k=y_{nk}+z_{nk}} \big\{ n_y,
n_z\big\}$$ where
\begin{eqnarray*}
n_y & = & \Big( \sum_{n \in \Z} \la^{-n\theta} \sum_{k=1}^\infty
\big\| a^{q/p_0}y_{nk}b^{q'/p_0'} \big\|_2^2 \Big)^\frac12, \\ n_z
& = & \Big( \sum_{n \in \Z} \la^{n(1-\theta)} \sum_{k=1}^\infty
\big\| a^{q/p_1}z_{nk}b^{q'/p_1'} \big\|_2^2 \Big)^\frac12.
\end{eqnarray*}
Using again the o.s.s. of
$\ell_2^{c_{p_0}}(w_\la^{-\theta};\ell_2)$ and
$\ell_2^{c_{p_1}}(w_\la^{1-\theta};\ell_2)$ as above, we get
$$\Big( \sum_{k=1}^\infty \|ax_kb\|_2^2 \Big)^\frac12 \lesssim
c_1(\la,\theta)^{-1} \, \big\| \mathbf{1}_\Z \otimes x
\big\|_{M_m(\mathcal{C}_K)}$$ where $\mathbf{1}_\Z = \sum_{n \in
\Z} \delta_n$ is the constant-$1$ sequence in $\Z$. Let us now
show that
$$\mathcal{C}_{c_{p_0'},c_{p_1'}}^J(w_\la^{-1},\theta) \subset
C_{q'}.$$ Indeed, arguing by homogeneity we assume that
$\sum_{n,k} x_{nk} \ten e_{(n,k),1}$ satisfies $$\Big\| \sum_{n,k}
x_{nk} \ten e_{(n,k),1} \Big\|_{M_m(\mathcal{C}_J)} < 1.$$ That
is, we may find $(z_{nk})$ such that $\sum_n x_{nk} = \sum_n
z_{nk}$ and $$\max \Big\{ \Big\| \sum_{n,k}
\la^{\frac{n\theta}{2}} z_{nk} \ten e_{(n,k),1}
\Big\|_{M_m(C_{p_0'})}, \Big\| \sum_{n,k}
\la^{\frac{-n(1-\theta)}{2}} z_{nk} \ten e_{(n,k),1}
\Big\|_{M_m(C_{p_1'})} \Big\} \le 1.$$ Taking $z_n = \sum_k z_{nk}
\otimes e_k$ and $$A = \mathcal{L}_{a^{q'/p_1'}}
\mathcal{R}_{b^{q/p_1}} \quad \mbox{and} \quad B =
\mathcal{L}_{a^{q'/p_0'}} \mathcal{R}_{b^{q/p_0}},$$ we observe
that
\begin{eqnarray*}
\lefteqn{\hskip-15pt \sum_{n \in \Z} \la^{-n(1-\theta)}
\|Az_n\|_2^2 + \sum_{n \in \Z} \la^{n\theta} \|Bz_n\|_2^2} \\
& = & \sum_{n,k} \la^{-n(1-\theta)} \big\| a^{q'/p_1'} z_{nk}
b^{q/p_1} \big\|_2^2 + \sum_{n,k} \la^{n\theta} \big\| a^{q'/p_0'}
z_{nk} b^{q/p_0} \big\|_2^2 \ \le \ 2
\end{eqnarray*}
holds by our assumption. According to Lemma \ref{dis1}, $\xi =
A^{2\theta-1}B^{1-2\theta} \summ_n z_n$ satisfies $\|A^{1-\theta}
B^{\theta} \xi\|_2 \lesssim c_2(\la,\theta)$ and for $x = \sum_k
x_k \ten e_{k,1}$ with $x_k = \sum_n x_{nk} = \sum_n z_{nk}$, we
have
$$\Big( \sum_{k=1}^\infty \|ax_kb\|_2^2 \Big)^\frac12 = \big\|
A^{\theta} B^{1-\theta}x \big\|_2 = \Big\| A^{1-\theta} B^{\theta}
A^{2\theta-1} B^{1-2\theta} \big( \sum_{n \in \Z} z_n \big)
\Big\|_2 \lesssim c_2(\la,\theta).$$ Therefore, duality yields
$C_q \subset \mathcal{C}_{c_{p_0},c_{p_1}}^K(w_\la,\theta)$ and
the assertion follows. \fin

Our next step is to construct a complete embedding of $S_q(X)$
into a 4-term sum. Together with Proposition \ref{Kvers} below,
this cb-embedding will be the key towards the main result in this
section. Let $\mathcal{K}_{p,q}(w;X)$ be defined by
$$S_p(X) \, + \, C_p \ten_h X \ten_h \ell_2^{r_q}(w) \, + \,
\ell_2^{c_q}(w) \ten_h X \ten_h R_p \, + \, \ell_2^{c_q}(w) \ten_h
X \ten_h \ell_2^{r_q}(w).$$ Let us write $\mathcal{K}_{p,q}(w)$
for the same space when $X=\C$ and $\mathcal{K}_{p,q}(w_\lambda;
X)$ / $\mathcal{K}_{p,q}(w_\lambda)$ for exponential sequences.
Here it is important to recall that we will be considering weights
$w$ on the index set $\Z \times \N$ which are constant on the
$\N$-component, so that (using the terminology above) another
description for this space could be
$$\mathcal{K}_{p,q}(w;X) = \Big[ C_p(\Z \times \N) +
\ell_2^{c_q}(w;\ell_2) \Big] \otimes_h X \otimes_h \Big[ R_p(\Z
\times \N) + \ell_2^{r_q}(w; \ell_2) \Big].$$ In the following
result, we study a map $S_q(X) \to \mathcal{K}_{p_0,p_1}(w;X)$ of
the form $$\sum_{k, \ell =1}^\infty e_{k,1} \otimes x_{k\ell}
\otimes e_{1,\ell} \mapsto \sum_{i,j=-\infty}^\infty
w_{ij}(p_0,p_1,q) \sum_{k,\ell=1}^\infty e_{i,1} \otimes e_{k,1}
\otimes x_{k\ell} \otimes e_{1,\ell} \otimes e_{1,j}.$$ Just to
shorten the notation, we change the order of tensors and write $$x
\mapsto \Big( \sum_{i,j=-\infty}^\infty w_{ij}(p_0,p_1,q) \,
e_{ij} \Big) \otimes x.$$ With this terminology, we have
$\mathbf{1}_\Z \otimes \mathbf{1}_\Z = \sum_{i,j \in \Z} e_{ij}$
for $\mathbf{1}_\Z = \sum_{n \in \Z} \delta_n$ as above.

\begin{proposition} \label{emb1}
If $p_0<q<p_1$ with $\frac{1}{q} = \frac{1-\theta}{p_0} +
\frac{\theta}{p_1}$ and $\la>1$, then $$u: x \in S_q(X) \mapsto
\Big( \sum_{i,j=-\infty}^\infty \la^{-(i+j)\theta/2} \, e_{ij}
\Big) \ten x \in \mathcal{K}_{p_0,p_1}(w_\la;X)$$ is a completely
isomorphic embedding with constants depending only on
$(\la,\theta)$.
\end{proposition}

\dem According to Lemma \ref{Xu}, the mappings
\begin{eqnarray*}
x \in C_q & \mapsto & \mathbf{1}_\Z \otimes x \in
\mathcal{C}_{c_{p_0},c_{p_1}}^K(w_\la,\theta), \\ x \in R_q &
\mapsto & \mathbf{1}_\Z \otimes x \in
\mathcal{C}_{r_{p_0},r_{p_1}}^K(w_\la,\theta),
\end{eqnarray*}
are cb-isomorphisms. Recalling that $\mathbf{1}_\Z \otimes
\mathbf{1}_\Z = \sum_{i,j \in \Z} e_{ij}$, we get $$x \in C_q
\ten_h X \ten_h R_q \mapsto \Big( \summ_{i,j}^{\null} e_{ij} \Big)
\otimes x \in \mathcal{C}_{c_{p_0},c_{p_1}}^K(w_\la,\theta) \ten_h
X \ten_h \mathcal{C}_{r_{p_0},r_{p_1}}^K(w_\la,\theta)$$ a
complete isomorphism. The right hand side inherits its o.s.s. from
$$\mathcal{G}_{c_{p_0},c_{p_1}}^K(w_\la,\theta) \ten_h X \ten_h
\mathcal{G}_{r_{p_0},r_{p_1}}^K(w_\la,\theta) \ = \ \sum_{i,j=1,2}
\mathcal{U}_i \ten_h X \ten_h \mathcal{V}_j$$ $$\mathcal{U}_1 =
\ell_2^{c_{p_0}}(w_\la^{-\theta};\ell_2) \ , \ \mathcal{U}_2 =
\ell_2^{c_{p_1}}(w_\la^{1-\theta};\ell_2) \ , \ \mathcal{V}_1 =
\ell_2^{r_{p_0}}(w_\la^{-\theta};\ell_2) \ , \ \mathcal{V}_2 =
\ell_2^{r_{p_1}}(w_\la^{1-\theta};\ell_2).$$

\noindent Thus, is suffices to show that the map $$z \in
\sum_{i,j=1,2} \mathcal{U}_i \ten_h X \ten_h \mathcal{V}_j \mapsto
\Big( \sum_{i \in \Z} \la^{-i\theta/2} e_{ii} \Big) \, z \, \Big(
\sum_{j \in \Z} \la^{-j\theta/2} e_{jj} \Big) \in
\mathcal{K}_{p_0,p_1}(w_\lambda;X)$$ is a complete embedding, in
this case with constants independent on $\la$ and $\theta$.
Moreover, since both spaces are the sum of $4$ spaces indexed
respectively by $(p_0,p_0)$, $(p_0,p_1)$, $(p_1,p_0)$ and
$(p_1,p_1)$, it clearly suffices to check our claim term by term.
However, this later fact follows from repeated use of the complete
isometries
\begin{eqnarray*}
\lefteqn{\hskip-30pt z \in \ell_2^{c_p}(w_{\la_1};\ell_2) \ten_h X
\ten_h \ell_2^{r_q}(w_{\la_2};\ell_2)} \\ [5pt] & \mapsto & \Big(
\sum_{i \in \Z} \la_1^{i/2} e_{ii} \Big) \, z \, \Big( \sum_{j \in
\Z} \la_2^{j/2} e_{jj} \Big) \in C_p \ten_h X \ten_h R_q,
\end{eqnarray*}
with $\lambda_1, \lambda_2 \in \{\lambda^{-\theta},
\lambda^{1-\theta}\}$ and $p,q \in \{p_0,p_1\}$. Details are left
to the reader. \fin

\begin{remark}
\emph{The cb-embedding of $L_p(\mathcal{M})$ into a von Neumann
algebra predual from \cite{JP3,JP4} can be described by means of
the map $u: L_p(\mathcal{M}) \to \mathcal{K}_{1,2}(w_\la)$ defined
on $\mathcal{M} \bar\ten \mathcal{B}(\ell_2(\Z))$ with the weight
given by $\la$. Indeed, it suffices to apply the Poisson map from
\cite[Section 8.2]{JP4} (a suitable average of sums of independent
copies which embeds $L_1+L_2^{r+c}(w)$ in $L_1$, for a suitable
strictly seminfinite weight $w$) with coefficients in $OH$ in
order to embed $\mathcal{K}_{1,2}(w_\la)$ into some
$L_1(\mathcal{A})$.}
\end{remark}

We shall also need an extended form of the embedding of
Proposition \ref{emb1} for arbitrary von Neumann algebras. More
concretely, that we have an isomorphic embedding
$L_{p_0}(\mathcal{M}; S_q) \to
L_{p_0}(\mathcal{M};\mathcal{K}_{p_0,p_1}(w_\lambda))$.
Fortunately, we only need this in the scalar case $X = \C$,
something that simplifies our approach very much. Our first task
is to define the space
$L_{p_0}(\mathcal{M};\mathcal{K}_{p_0,p_1}(w))$ appropriately. We
have
\begin{eqnarray*}
\mathcal{K}_{p_0,p_1}(w) & = & \Big[ C_{p_0} + C_{p_1}(w) \Big]
\otimes_h \Big[ R_{p_0} + R_{p_1}(w) \Big],
\end{eqnarray*}
where the row/column spaces are taken in the index set $\Z \times
\N$ and the weight $w$ is constant on the $\N$-component. Recall
that in Section \ref{S1} we have defined the spaces
$L_{2p}^r(\mathcal{M}; C_q)$ and $L_{2p}^c(\mathcal{M}; R_q)$ and
the same definition is valid for the weighted row/column spaces.
Moreover, we may define
\begin{eqnarray*}
L_{2p}^r(\mathcal{M}; C_{q_1}(w_1) + C_{q_2}(w_2)) & = &
L_{2p}^r(\mathcal{M}; C_{q_1}(w_1)) + L_{2p}^r(\mathcal{M};
C_{q_2}(w_2)), \\ L_{2p}^c(\mathcal{M}; R_{q_1}(w_1) +
R_{q_2}(w_2)) & = & L_{2p}^c(\mathcal{M}; R_{q_1}(w_1)) +
L_{2p}^c(\mathcal{M}; R_{q_1}(w_2)),
\end{eqnarray*}
for arbitrary weights by viewing them as both embedded into the
space of sequences with values in $L_{2p}(\mathcal{M})$, and
thereby defining the sum by taking the corresponding quotients.
This allows us to consider
$$L_{p_0}(\mathcal{M};\mathcal{K}_{p_0,p_1}(w)) =
L_{2p_0}^r \big( \mathcal{M}; C_{p_0} + C_{p_1}(w) \big)
\ten_{\mathcal{M};h} L_{2p_0}^c \big( \mathcal{M}; R_{p_0} +
R_{p_1}(w) \big)$$ with norm given by
\begin{eqnarray*}
\inf_{x_{ij} = \sum_k \al_{ik}\beta_{kj}} \big\| (\alpha_{ik} )
\big\|_{L_{2p_0}^r(\mathcal{M}; C_{p_0} + C_{p_1}(w)) \ten_h R}
\big\| ( \beta_{kj} ) \big\|_{C \ten_h L_{2p_0}^c(\mathcal{M};
R_{p_0} + R_{p_1}(w))}.
\end{eqnarray*}

\begin{proposition} \label{emb1C}
If $p_0<q<p_1$ with $\frac{1}{q} = \frac{1-\theta}{p_0} +
\frac{\theta}{p_1}$ and $\la>1$, then $$u: x \in
L_{p_0}(\mathcal{M};S_q) \mapsto \Big( \sum_{i,j=-\infty}^\infty
\la^{-(i+j)\theta/2} \, e_{ij} \Big) \ten x \in
L_{p_0}(\mathcal{M}; \mathcal{K}_{p_0,p_1}(w_\la))$$ is a complete
isomorphic embedding with constants depending only on
$(\la,\theta)$.
\end{proposition}

\dem Lemma \ref{Xu} remains valid here as well, i.e. we have
\begin{eqnarray*}
L_{2p_0}^r(\mathcal{M}; \mathcal{C}_{c_{p_0},
c_{p_1}}^K(w_\la,\theta)) & \simeq_{cb} & L_{2p_0}^r(\mathcal{M}; C_q), \\
L_{2p_0}^c(\mathcal{M}; \mathcal{C}_{r_{p_0},
r_{p_1}}^K(w_\la,\theta)) & \simeq_{cb} & L_{2p_0}^c(\mathcal{M};
R_q).
\end{eqnarray*}
Indeed, by the factorization properties of the spaces involved it
really suffices to prove this for $p_0=\infty$, and then the exact
same argument in Lemma \ref{Xu} applies since the key formula is
the operator space structure of $C_q$, which according to
\cite{H4,P5} is still valid for arbitrary von Neumann algebras
$$\Big\| \summ_k x_k \ten e_{k,1} \Big\|_{L_{\infty}(\mathcal{M}; C_q)}
= \sup_{\|a\|_{L_{2q}(\mathcal{M})}, \|b\|_{L_{2q'}(\mathcal{M})}
\le 1} \Big( \summ_k \|a x_k b\|_{L_2(\mathcal{M})}^2
\Big)^\frac12.$$ Since we have $$L_{p_0}(\mathcal{M}; S_q) =
L_{2p_0}^r(\mathcal{M}; C_q) \ten_{\mathcal{M},h}
L_{2p_0}^c(\mathcal{M}; R_q),$$ an element in
$L_{p_0}(\mathcal{M}; S_q)$ factorizes as a product of two
elements in $L_{2p_0}^r(\mathcal{M}; C_q)$ and
$L_{2p_0}^c(\mathcal{M}; R_q)$ respectively. Thus, we deduce it
can be written as a product of two sums from
$$L_{2p_0}^r(\mathcal{M}; \mathcal{C}_{c_{p_0},
c_{p_1}}^K(w_\la,\theta)) \quad \mbox{and} \quad
L_{2p_0}^c(\mathcal{M}; \mathcal{C}_{r_{p_0},
r_{p_1}}^K(w_\la,\theta))$$ respectively. Arguing as in
Proposition \ref{emb1}, we see that $u$ is bounded. To prove the
converse, we observe that a norm estimate for $u(x)$ means a
factorization of the form $x_{k\ell} = \sum_m \la^{i\theta/2}
\al_{ik,m} \beta_{m,j\ell} \la^{j\theta/2}$ valid for all $i,j \in
\Z$ and with
\begin{eqnarray*}
(\al_{ik,m}) & \in & L_{2p_0}^r(\mathcal{M}; C_{p_0} +
C_{p_1}(w_\lambda)) \ten_h R, \\ (\beta_{m,j \ell}) & \in & C
\ten_h L_{2p_0}^c(\mathcal{M}; R_{p_0} + R_{p_1}(w_\lambda)).
\end{eqnarray*}
We may rewrite this as $$x_{k\ell} = a_{ik} b_{j\ell}$$ where
$a_{ik}=\la^{i\theta/2}\sum_m \al_{ik,m} \ten e_{1m}$ and
$b_{j\ell} = \la^{j\theta/2} \sum_m \beta_{j\ell} \ten e_{m1}$. We
may also assume by approximation that we are only dealing with
finitely many nonzero entries $x_{k\ell}$ with full left and right
support in a finite von Neumann algebra $\mathcal{M}$. Let
$e_{j\ell}$ be the left support of $b_{j\ell}$. Then we deduce
from $a_{ik}b_{j\ell} = x_{k\ell} = a_{i'k}b_{j\ell}$ that we have
$a_{ik}e_{jl}= a_{i'k}e_{jl}$. This holds for all indices $j,\ell$
and hence for $e = \vee  e_{j\ell}$ we find $a_{ik}e= a_{i'k}e$
for all $k$ and $i\neq i'$. Taking $v_k=a_{ik}e$, we deduce
$$x_{k\ell} = v_k b_{j\ell}$$ for all $j$. Similarly, let $f$ be
supremum of the right supports of the $v_k$'s. As above we deduce
that $fb_{j\ell}=fb_{j'\ell}$ for all $\ell$ and $j \neq j'$. Thus
we may define $w_\ell = fb_{j\ell}$ and obtain a factorization
$$x_{k\ell} = v_kw_\ell$$ such that $v_k= a_{ik}e$ and $w_\ell =
fb_{j\ell}$. Since the space $L_{2p_0}^r(\mathcal{R}; C_{p_0} +
C_{p_1}(w_\lambda))$ is a right $\mathcal{R}$-module and
$L_{2p_0}^c(\mathcal{R}; R_{p_0} + R_{p_1}(w_\lambda))$ is a left
$\mathcal{R}$-module, we may now apply the announced extension of
Lemma \ref{Xu} and deduce that
$$(v_k) \in L_{2p_0}^r(\mathcal{M}; C_q) \quad , \quad (w_\ell) \in
L_{2p_0}^c(\mathcal{M};R_q).$$ This implies $x=(x_{kl}) \in
L_{p_0}(\mathcal{M}; S_q)$ and hence $u$ is an isomorphism.
Tensoring with another copy of $L_{p_0}(M_n)$ does not change
constants in this argument and hence $u$ is indeed a complete
isomorphism. The proof is complete. \fin

\subsection{Change of density}

We need an alternative description of $L_p(\mathcal{M};
\mathcal{K}_{q_1,q_2}(w))$ according to another description of the
space $L_p(\mathcal{M}; C_{q_1} \ten_h R_{q_2})$. Namely if we
take $\frac1p = \frac{1}{s_j} + \frac{1}{q_j}$ and $\frac{1}{q} =
\frac{1}{2q_1} + \frac{1}{2q_2}$, we have the Banach space
isometry $$L_p(\mathcal{M}; C_{q_1} \ten_h R_{q_2}) =
L_{\underline{2s_1}q\underline{2s_2}}(\mathcal{R}; \mathcal{M}).$$
This follows again by complex interpolation. In particular,
$L_p(\mathcal{M}; \mathcal{K}_{q_1,q_2}(w))$ is Banach space
isomorphic to a 4-term sum of amalgamated $L_p$ spaces. Certain
embedding in \cite{JP4} for these spaces will be essential in the
following change of density argument. Recall the notion of
$(p_1,p_2)$-convex maps $T: X \to L_p(\mathcal{M})$ from the
Introduction.

\begin{proposition} \label{Kvers}
Let $1 \le p < p_1 \wedge p_2 \le \infty$ and $$\al = \Big(
\frac{1}{p}-\frac{1}{p_2} \Big) \Big/ \Big(
\frac{1}{p}-\frac{1}{p_1} \Big).$$ If $T: X \to L_p(\mathcal{M})$
is $(p_1,p_2)$-convex and $w$ is any weight, then $$T \ten id:
\mathcal{K}_{p,p_1}(w;X) \to L_p \big( \mathcal{M};
\mathcal{K}_{p,p_2}(w^{\al}) \big) \quad \mbox{with} \quad
(w^\alpha)_n = (w_n)^\alpha$$ is completely bounded and its
cb-norm can be estimated by $c(p,p_2) \, k_{(p_1,p_2)}(T)$.
\end{proposition}

\dem We may and will assume that $X$ is finite dimensional. Given
natural numbers $m,n \in \N$, consider a faithful state $\phi_m$
on $M_m$ with density $d_{\phi_m}$. We may regard the density
$d_{n \phi_m}$ of $n \phi_m$ as a diagonal operator whose entries
form a weight on the index set $\{1,2,\ldots,m\}$. Define
$$\mathcal{K}_{p,q}^n(\phi_m;X) = \Big[
\ell_2^{c_p}(d_{n\phi_m}^\frac1p) + \ell_2^{c_q}(d_{n
\phi_m}^\frac1q) \Big] \ten_h X \ten_h \Big[
\ell_2^{r_p}(d_{n\phi_m}^\frac1p) + \ell_2^{r_q}(d_{n
\phi_m}^\frac1q) \Big].$$ As above, the expression
$\mathcal{K}_{p,q}^n(\phi_m)$ will be reserved for the
scalar-valued case. The space $\mathcal{K}_{p,q}^n(\phi_m;X)$ can
be written as a 4-term sum of asymmetric $L_p$ spaces as in
\cite{JP1,JP3,JP4}. Namely, if we consider the asymmetric spaces
$$L_{(2p,2q)}(\phi_m;X) = \ell_2^{c_p}(d_{\phi_m}^\frac1p) \ten_h
X \ten_h \ell_2^{r_q}(d_{\phi_m}^\frac1q),$$ we have $$\Big\|
\sum_{i,j=1}^m x_{ij} \ten e_{ij} \Big\|_{L_{(2p,2q)}(n\phi_m;X)}
= n^{\frac{1}{2p}+\frac{1}{2q}} \Big\| d_{\phi_m}^{\frac{1}{2p}}
\Big( \sum_{i,j=1}^m x_{ij} \ten e_{ij} \Big)
d_{\phi_m}^{\frac{1}{2q}} \Big\|_{C_p^m \ten_h X \ten_h R_q^m}.$$
This gives a description of $\mathcal{K}_{p,q}^n(\phi_m;X)$ in
terms of asymmetric Schatten classes.

\noindent Consider the $n$-fold free product
$$\mathcal{A}_n = (M_m,\phi_m)^{*n}.$$ According to \cite{J2}, we
know that $\mathcal{A}_n$ is QWEP. In particular, it is very
well-known the existence of a normal $*$-homomorphism $\rho$ and a
normal conditional expectation $\mathcal{E}$ as follows $$\rho:
\mathcal{A}_n \to \Big( \prodd_\U S_1 \Big)^* \quad \mbox{and}
\quad \E: \Big( \prodd_\U S_1 \Big)^* \to \mathcal{A}_n.$$ We also
know that we have $L_p$ extensions $\rho_p$ and $\E_p$ for $1 \le
p < \infty$. Let us denote by $\pi_j: M_m \to \mathcal{A}_n$ the
$j$-th coordinate map. Then $\rho \hskip1pt \pi_j: M_m \to [
\hskip1pt \prod_\U S_1 \hskip1pt ]^*$ is a $*$-homomorphism.
Following an argument of Kirchberg, we observe that by Kaplansky's
density theorem the unit ball of $\prod_\U S_{\infty}$ is dense in
the strong and strong$^*$ topology of $[ \hskip1pt \prod_\U
S_1]^*$. Let $B=\ell_{\infty}^{st*}({\mathcal I},\prod
S_{\infty})$ the $C^*$-algebra of all strong and strong$^*$
converging families. Then $[ \hskip1pt \prod_\U S_1]^*$ is a
quotient of $B$. Since $M_m$ is nuclear, we can apply the
Choi-Effros lifting theorem \cite[Theorem 3.10]{CE} for the maps
$\pi_j$ and find nets $v_{s,j}: M_m \to \prod_\U S_{\infty}$ of
completely positive and contractive maps such that $(v_{s,j})$
converges to $\rho \pi_j$ in the strong and strong$^*$ topologies.
Let us consider the maps $$u_1: x \in
\mathcal{K}_{p,p_1}^n(\phi_m;X) \mapsto \sum_{j=1}^n \rho_p
\pi_j(x) \ten \delta_j \in \prodd_\U S_p (\ell_{p_1}^n(X)),$$
$$u_2: x \in L_p(\mathcal{M}; \mathcal{K}_{p,p_2}^n(\phi_m))
\mapsto \sum_{j=1}^n (id_{L_p(\mathcal{M})} \otimes \pi_j)(x) \ten
\delta_j \in L_p(\mathcal{A}_n \bar\ten \mathcal{M};
\ell_{p_2}^n).$$ We claim that $u_1$ is completely contractive and
$u_2$ is an embedding. Let us note that $u_1$ is also a
cb-embedding, a fact which will not be needed nor proved in this
paper. The proof that $u_2$ is an embedding (defining
$L_p(\mathcal{M}; \mathcal{K}_{p,p_2}^n(\phi_m))$ as indicated
before the statement of this result) was given in Theorem 7.3 and
Remark 7.4 of \cite{JP4}. Moreover, we know that $u_2$ is a
complete contraction while the cb-norm of its inverse is
controlled by a constant $c(p,p_2)$, see Remark \ref{Const} below
for more on the value of $c(p,p_2)$. For the first part of the
claim, let us show that
$$\big\| u_1(x) \big\|_{\prod_\U S_p(\ell_{p_1}^n(X))} \le
n^\frac1p \|x\|_{L_{(2p,2p)}(\phi_m;X)}.$$ In fact, we will only
prove this inequality since the remaining ones for the terms
associated to $(2p,2p_1),(2p_1,2p)$ and $(2p_1,2p_1)$ are similar.
Indeed, we refer the reader to \cite[Proposition 3.5]{JP1} for the
exact same argument. Since $p < p_1$, we have
$$\big\| u(x) \big\|_{\prod_\U S_p(\ell_{p_1}^n(X))} \le \Big(
\sum_{j=1}^n \big\| \rho_p \pi_j(x) \big\|_{\prod_\U S_p(X)}^p
\Big)^\frac1p.$$ Therefore, it suffices to consider a fixed
component $j$. We may write $x=ayb^*$ such that $a,b \in
L_{2p}^r(\phi_m)$ are of norm $1$ and
$\|x\|_{L_{(2p,2p)}(\phi_m;X)} \sim \|y\|_{M_m(X)}$. Then, the
element $y_{s,j}$ defined by
$$y_{s,j} = \big( v_{s,j} \ten id_X \big)(y) \in \prodd_\U
S_{\infty}(X)$$ satisfies $\|y_{s,j}\|\le \|y\|_{M_m(X)}$.
Moreover, the strong convergence guarantees the norm convergence
of $\limm_s \, \rho_{2p}\pi_j(a) \, y_{s,j} \, \rho_{2p}\pi_j(b) =
\rho_{2p} \pi_j(a) \, \rho \pi_j(y) \, \rho_{2p} \pi_j(b) = \rho_p
\pi_j(x)$ (see \cite{JS} for further details) and we obtain
$$\big\| \rho_p \pi_j(x) \big\|_{\prod_\U S_p(X)} \le
\|x\|_{L_{(2p,2p)}(\phi_m;X)}.$$ Since the same inequality holds
after tensorizing with the identity on $S_p$, this proves our
claim. On the other hand, using the $(p_1,p_2)$-convexity of $T$
in conjunction with the contractivity of $u_1$, we deduce
\begin{eqnarray*}
\lefteqn{\hskip-6cm \Big\| \sum_{j=1}^n \rho_p \Big(
\underbrace{\pi_j \ten id_X \big( d_{n \phi_m}^{1/2p}(id_{M_m}
\ten T(x)) d_{n \phi_m}^{1/2p} \big)}_{\pi_j(Tx) \ \mathrm{for \
short}} \Big) \ \ten \delta_j \Big\|_{\prod_{\U}
S_p(L_p(\mathcal{M}; \ell_{p_2}^n))}} \\ \null \hskip6cm \null &
\le & k_{(p_1,p_2)}(T) \, \|x\|_{\mathcal{K}_{p,p_1}^n(\phi_m;X)}.
\end{eqnarray*}
Moreover, we may understand this as a cb-inequality, which remains
true after tensorizing with $id_{S_p}$. Then we recall from
\cite[Chapter 3]{JP4} that the space $L_p(\ell_{p_2})$ is stable
under the conditional expectation
$$\mathcal{E}_p: \prodd_\U S_p(L_p(\mathcal{M}; \ell_{p_2}^n))
\to L_p(\mathcal{A}_n \bar\ten \mathcal{M}; \ell_{p_2}^n).$$
Therefore, we have proved that $$\big\| \E_p u_1T:
\mathcal{K}_{p,p_1}^n(\phi_m;X) \to L_p(\mathcal{A}_n \bar\ten
\mathcal{M}; \ell_{p_2}^n) \big\|_{cb} \le k_{(p_1,p_2)}(T).$$
Note that the range of $\E_p u_1T$ is still of the form $$\E_p
u_1T(x)=\sum_{j=1}^n \pi_j(Tx)\ten \delta_j.$$ This means in
particular that $\E_p u_1T$ maps $\mathcal{K}_{p,p_1}^n(\phi_m;X)$
in the range of $$u_2 \hskip1pt [ \hskip1pt L_p(\mathcal{M};
\mathcal{K}_{p,p_2}^n(\phi_m)) \hskip1pt ].$$ Thus we obtain
$$\big\| T \ten id: \mathcal{K}_{p,p_1}^n(\phi_m;X) \to
L_p(\mathcal{M}; \mathcal{K}_{p,p_2}^n(\phi_m)) \big\|_{cb}
\lesssim c(p,p_2) \, k_{(p_1,p_2)}(T).$$ Let us now prove the
assertion. First we may replace $\phi_m$ by the state $\phi_m \ten
\tau_\ell$ on $M_{m\ell}$ where $\tau_\ell$ is the normalized
trace on $M_\ell$. Then we note that the space of elements $x \ten
e$, with $e$ a fixed projection satisfying $\tau_\ell(e) =
\gamma$, is simultaneously complemented in all the asymmetric
spaces $L_{(2p,2q)}$ considered. Thus, we restrict our attention
to this subspace. Moreover, we clearly have
$$n^{\frac{1}{2p}+\frac{1}{2q}} \|x\ten e\|_{L_{(2p,2q)}(\phi_m
\ten \tau_\ell;X)} = \big\| d_{n \gamma \phi_m}^{\frac{1}{2p}} x
\hskip1pt d_{n \gamma \phi_m}^{\frac{1}{2q}} \big\|_{C_p^m \ten_h
X \ten_h R_q^m}.$$ By \cite[Lemma 1.2]{JP1}, tensorizing with
$id_{S_{(2p,2q)}}$ we obtain a complete isometry $$x \ten e \in
\mathcal{K}_{p,p_1}^n(\phi_m \ten \tau_\ell;X) \mapsto d_{n \gamma
\phi_m}^{\frac{1}{2p}} x \hskip1pt d_{n \gamma
\phi_m}^{\frac{1}{2p}} \in \mathcal{K}_{p,p_1}(w_\lambda;X)$$ with
$w_\lambda = (n \gamma d_{\phi_m})^{\frac{1}{p_1}-\frac{1}{p}}$. A
similar argument leads to the complete isometry $$Tx \ten e \in
L_p \big( \mathcal{M}; \mathcal{K}_{p,p_2}^n(\phi_m \ten
\tau_\ell) \big) \mapsto d_{n \gamma \phi_m}^{\frac{1}{2p}} T
\hskip-1pt x \hskip1pt d_{n \gamma \phi_m}^{\frac{1}{2p}} \in L_p
\big( \mathcal{M}; \mathcal{K}_{p,p_2}(w_\mu) \big)$$ with $\mu =
(n \gamma d_{\phi_m})^{\frac{1}{p_2}-\frac{1}{p}} = w^\alpha$.
This implies the assertion for $w=(n \gamma
\phi_m)^{\frac{1}{p_1}-\frac{1}{p}}$. It just remains to show that
the general case follows from this one. Indeed, by approximation
it clearly suffices to show it for $w$ being a weight on
$\{1,2,\ldots,m\}$ as far as we see that the constants are
independent of $m$. Therefore, we have to see that every $w$
supported on $\{1,2,\ldots,m\}$ can be obtained in this form.
Given such a weight $w$, we consider the functional on $M_m$ given
by $$\psi_m \Big( \sum_{i,j=1}^m \alpha_{ij} e_{ij} \Big) =
\sum_{k=1}^m w_k^{\frac{pp_1}{p-p_1}} \alpha_{kk}$$ and the state
$\phi_m$ defined by $\psi_m = \psi_m(\mathbf{1}_{M_m}) \phi_m$.
Let us set $n = [\psi_m(\mathbf{1}_{M_m})] +1$ where $[ \cdot ]$
stands for the integer part. Let $0 < \gamma < 1$ be determined by
the relation $n \gamma = \psi_m(\mathbf{1}_{M_m})$. We may assume
by approximation that $\gamma$ is a rational number. Let
$\tau_\ell$ be the normalized trace on $M_\ell$. Taking $\ell$
large enough, we may consider a projection $e$ in $M_\ell$
satisfying $\tau_\ell(e)=\gamma$. Hence, the embedding $$x \in
(M_m,\phi_m) \mapsto x \ten e \in (M_{m \ell}, \phi_m \ten
\tau_\ell)$$ produces the desired identification $w = (n \gamma
d_{\phi_m})^{\frac{1}{p_1}- \frac1p}$. The proof is complete. \fin

\begin{remark}
\emph{The embedding of the algebra $\mathcal{A}_n$ into an
ultraproduct algebra is the key tool used in \cite{J4} to
generalize vector-valued noncommutative $L_p$ spaces to QWEP
algebras and such notion underlies the proof of Proposition
\ref{Kvers}. Note however that we do not need at any rate to
require $\mathcal{M}$ to be QWEP.}
\end{remark}

\begin{remark} \label{Const}
\emph{According to Remarks 2.2 and 5.7 of \cite{JP5}, the value of
the constant $c(p,p_2)$ above remains uniformly bounded in $p$ and
$p_2$ as far as $(p,p_2) \nsim (1, \infty)$. In that case, we only
know that it is controlled by $1+\frac{p_2-p}{pp_2+p-p_2}$. Note
that this singularity near $(1,\infty)$ seems to be removable
since the corresponding complete embedding holds at the point
$(1,\infty)$.}
\end{remark}

\begin{remark}
\emph{Although not needed for our purposes in this paper, let us
point a generalization of Proposition \ref{Kvers} for potential
applications. Only for this remark we shall write
$L_p[\mathcal{M}; X]$ to denote the generalization of
$L_p(\mathcal{M}; X)$ for QWEP von Neumann algebras in \cite{J4}.
Assume that $1 \le s \le u \wedge v \le u \vee v < p_1 \wedge p_2
\le \infty$ and set $\beta = ( \frac{1}{s} - \frac{1}{p_2} ) / (
\frac{1}{u} - \frac{1}{p_1} )$. If the map $T: X \to
L_v(\mathcal{M})$ is $(p_1,p_2)$-convex and $w$ is any weight,
then $$T \ten id: \mathcal{K}_{u,p_1}(w;X) \to L_v \big(
\mathcal{M}; \mathcal{K}_{s,p_2}(w^\beta) \big)$$ is completely
bounded and its cb-norm can be estimated by $c(s,p_2) \,
k_{(p_1,p_2)}(T)$. The proof follows the same pattern. Indeed,
arguing as above we know that the mapping $\E_u T u_1:
\mathcal{K}_{u,p_1}^n(\phi_m;X) \to L_u[\mathcal{A}_n;
L_v(\mathcal{M}; \ell_{p_2}^n)]$ is completely bounded. Moreover,
we also have complete contractions
$$L_u \big[ \mathcal{A}_n; L_v(\mathcal{M}; \ell_{p_2}^n) \big]
\to L_s \big[ \mathcal{A}_n; L_v(\mathcal{M}; \ell_{p_2}^n) \big]
\to L_v \big[ \mathcal{M}; L_s(\mathcal{A}_n; \ell_{p_2}^n)
\big]$$ given by the identity map. The first one follows from the
fact that $s \le u$ and $\mathcal{A}_n$ is a noncommutative
probability space. The second one follows from Minkowski's
inequality since $s \le v$. Then, we use again the embedding
$$L_v \big( \mathcal{M};\mathcal{K}_{s,p_2}^n(\phi_m) \big) \to
L_v \big[ \mathcal{M}; L_s(\mathcal{A}_n; \ell_{p_2}^n) \big]$$ to
conclude
$$\big\| T \ten id: \mathcal{K}_{u,p_1}^n(\phi_m;X) \to L_v \big(
\mathcal{M}; \mathcal{K}_{s,p_2}^n(\phi_m) \big) \big\|_{cb}
\lesssim c(s,p_2) \, k_{(p_1,p_2)}(T).$$ The change of density in
this case is given by $$w = (n \gamma d_{\phi_m})^{\frac{1}{p_1} -
\frac{1}{u}} \quad \mbox{and} \quad \mu = (n \gamma
d_{\phi_m})^{\frac{1}{p_2} - \frac{1}{s}}.$$ Thus, it turns out
that $\mu = w^\beta$ for our choice of $\beta$. This completes the
argument.}
\end{remark}

\noindent Now we are ready for the key embedding of this paper.

\vskip3pt

\demB Let $$\frac{1}{q} = \frac{1-\theta}{p} + \frac{\theta}{p_1}
= \frac{1-\eta}{p} + \frac{\eta}{p_2}$$ and $\la>1$. Then we have
the identity $$\al \eta = \frac{\frac{1}{p} -
\frac{1}{p_2}}{\frac{1}{p} - \frac{1}{p_1}} \, \frac{\frac{1}{p} -
\frac{1}{q}}{\frac{1}{p} - \frac{1}{p_2}} = \frac{\frac{1}{p} -
\frac{1}{q}}{\frac{1}{p} - \frac{1}{p_1}} = \theta,$$ where $\al$
is the real number defined in Proposition \ref{Kvers}. Let
$$u_{\theta,\la}: S_q(X) \to \mathcal{K}_{p,p_1}(w_\la;X) \quad
\mbox{and} \quad u_{\eta,\mu}: S_q \to
\mathcal{K}_{p,p_2}(w_\mu)$$ be the cb-embeddings given by
Proposition \ref{emb1}. Taking $\mu=\la^\al$, we note that $$(T
\ten id) u_{\theta,\la} = (u_{\eta,\mu} \ten
id_{L_p(\mathcal{M})}) (T \ten id).$$ Indeed, we deduce from
$\mu^{\eta}=\la^{\al \eta}=\la^{\theta}$ that
\begin{eqnarray*}
\lefteqn{\hskip-40pt (u_{\eta,\mu} \ten id_{L_p(\mathcal{M})}) (T
\ten id) (x) \ = \ \Big( \sum_{i,j=-\infty}^\infty
\mu^{-(i+j)\eta/2} \, e_{ij} \Big) \ten T(x)} \\ & = & \Big(
\sum_{i,j=-\infty}^\infty \la^{-(i+j)\theta/2} \, e_{ij} \Big)
\ten T(x) \ = \ (T \ten id) u_{\theta,\la}(x).
\end{eqnarray*}
According to Proposition \ref{Kvers}, we know that $$T \ten id:
\mathcal{K}_{p,p_1}(w_\la;X) \to L_p \big( \mathcal{M};
\mathcal{K}_{p,p_2}(w_\la^{\al}) \big)$$ is completely bounded and
hence $(T \ten id) u_{\theta,\la}$ is completely bounded. Thus, we
derive that $(u_{\eta,\mu} \ten id_{L_p(\mathcal{M})})(T \ten id)$
is completely bounded. Since $u_{\eta,\mu} \ten
id_{L_p(\mathcal{M})}$ is a cb-embedding by Proposition
\ref{emb1C}, we obtain \vskip-10pt \null \hfill $\qquad \big\| T
\ten id: S_q(X) \to L_p(\mathcal{M}; S_q) \big\|_{cb} \le
c(p,q,p_1,p_2) \, k_{(p_1,p_2)}(X)$. \hfill \fin

\vskip3pt

\begin{remark} \label{Constante}
\emph{Keeping track of constants, we have}
\begin{eqnarray*}
c(p,q,p_1,p_2) & \lesssim & \frac{pp_2}{pp_2 + p - p_2} \
\inf_{\lambda > 1} \frac{c_2(\lambda,\theta)^2}{c_1(\lambda^\alpha,
\theta/\alpha)^2}
\\ & = & \frac{pp_2}{pp_2 + p - p_2} \ \inf_{\lambda > 1}
\frac{(2\lambda - \lambda^\theta -
\lambda^{1-\theta})(\lambda^{\alpha
-\theta}-1)}{(\lambda^{1-\theta}-1)(\lambda^\theta +
\lambda^{\alpha-\theta}-2)} \\ & \le & \frac{pp_2}{pp_2 + p - p_2}
\ \lim_{\lambda \to 1^+} \frac{2\lambda - \lambda^\theta -
\lambda^{1-\theta}}{\lambda^{1-\theta}-1} \ = \ \frac{pp_2}{pp_2 +
p - p_2} \, \frac{1}{1-\theta},
\end{eqnarray*}
\emph{unless $(p,p_2) = (1,\infty)$, in which case the first term
on the right behaves like $1$.}
\end{remark}

\section{Maurey's factorization and applications}

We now prove an operator space form of Maurey's factorization
theorem. Then we will establish some selected applications for
operator spaces, noncommutative $L_p$ spaces and Fourier
multipliers.

\subsection{Maurey's factorization}

Let us begin with some basic inequalities to be used below. We
refer the reader to the Introduction for the definition of the
operator space analogs of cotype $p$ and absolutely
$(p,1)$-summing maps.

\begin{lemma} \label{prel}
Let $2 \le p \le \infty:$
\begin{enumerate}
\item[i)] If $T$ has cb-cotype p, then $$\pi_{p,1}^{cb}(T) \le
c_p^{cb}(T).$$

\item[ii)] $id_{L_p(\mathcal{M})}$ is completely $(p,1)$-summing
for any algebra $\mathcal{M}$.

\item[iii)] Let us consider two von Neumann algebras $\mathcal{M},
\mathcal{N}$ and assume that the map $T: L_q(\mathcal{M}) \to
L_p(\mathcal{N})$ is a completely bounded map. Then, the following
inequality holds for $1 \le q \le \infty$ $$\big\| T \ten id:
L_q(\mathcal{M}; \ell_1) \to \ell_p(L_p(\mathcal{N})) \big\|_{cb}
\le \|T\|_{cb}.$$
\end{enumerate}
\end{lemma}

\dem Consider $\Omega = \mathbb{T}^{\N}$ equipped with the product
topology and the corresponding Haar measure $\mu$. Clearly, the
map $j: \ell_1 \to C(\Omega)$ given by $j(\al)(\omega) = \sum_k
\omega_k \al_k$ is a complete contraction. Hence, we have
$$\big\| j \ten id_X: \ell_1 \ten_{\min} X \to
L_{\infty}(\Omega) \ten_{\min} X \big\|_{cb} \le 1.$$ The
inclusion $L_{\infty}(\Omega; X) \subset L_p(\Omega; X)$ is also
completely contractive and $$j \ten id_X(\ell_1 \ten_{\min} X)
\subset \Rad_p(X).$$ Hence i) follows by definition. To prove ii)
it suffices to show that the space $L_p(\mathcal{M})$ has
cb-cotype $p$. Let $\Lambda: f \in L_{\infty}(\Omega) \bar{\ten}
\mathcal{M} \mapsto (\int_\Omega f \, \eps_k \, d\mu)_{k \ge 1}
\in \ell_{\infty}(\mathcal{M})$ be the Rademacher coefficient map.
$\Lambda$ is a complete contraction and coincides with the
orthogonal projection $\Lambda: L_2(\Omega; L_2(\mathcal{M})) \to
\ell_2(L_2(\mathcal{M}))$. Thus, by interpolation we deduce that
$\Lambda: L_p(\Omega; L_p(\mathcal{M})) \to
\ell_p(L_p(\mathcal{M}))$ is a contraction. We conclude by
restriction to $\Rad_p(L_p(\mathcal{M}))$. Assertion iii) now
follows from the fact that the inclusion $L_q(\mathcal{M}; \ell_1)
\subset L_q(\mathcal{M}) \ten_{\min} C(\Omega)$ is completely
contractive. Indeed, in that case, we may compose with
\begin{eqnarray*}
L_q(\mathcal{M}) \ten_{\min} C(\Omega) &
\stackrel{T}{\longrightarrow} & L_p(\mathcal{N}) \ten_{\min}
C(\Omega) \\ & \stackrel{id}{\longrightarrow} &
L_p(\Omega;L_p(\mathcal{N})) \\ &
\stackrel{\Lambda}{\longrightarrow} & \ell_p(L_p(\mathcal{N})).
\end{eqnarray*}
It therefore suffices to show that for every $\omega \in \Omega$,
the map $$\phi_{\omega}: L_q(\mathcal{M}; \ell_1) \to
L_q(\mathcal{M}) \quad \mbox{with} \quad \phi_{\omega}(x) =
\summ_k \omega_k x_k$$ is completely contractive. Recall that
$S_q^m(L_q(\mathcal{M}; \ell_1)) = L_q(M_m \ten \mathcal{M};
\ell_1)$ and hence we just need to show that $\phi_w$ is a
contraction for all $w \in \Omega$. Assume $x_k=\sum_{j}
a_{kj}b_{kj}$ such that
$$\Big\| \big( \summ_{k,j} a_{kj} a_{kj}^* \big)^\frac12
\Big\|_{2q} \Big\| \big( \summ_{k,j} b_{kj}^*b_{kj} \big)^\frac12
\Big\|_{2q} \le 1.$$ Then, the Cauchy-Schwartz inequality implies
\\ [8pt] \hfill $\begin{array}{rcl} \displaystyle \Big\| \summ_k \omega_k
x_k \Big\|_q & = & \displaystyle \Big\| \summ_{k,j} \omega_k
a_{kj} b_{kj} \Big\|_q \\ [10pt] & \le & \displaystyle \Big\|
\big( \summ_{k,j} a_{kj} a_{kj}^* \big)^\frac12 \Big\|_{2q} \Big\|
\big( \summ_{k,j} |\omega_k|^2 b_{kj}^*b_{kj} \big)^\frac12
\Big\|_{2q} \le 1.
\end{array}$ \hfill $\begin{array}{c} \null \\ [13pt] \square \end{array}$

\vskip12pt

\begin{lemma} \label{elm}
Let $1\le p<\infty$. Then $L_p(\mathcal{M})$ has cb-cotype
$q=\max\{p,p'\}$.
\end{lemma}

\dem Let $\Om=\{-1,1\}^n$ with Haar measure $\mu$. Given $2 \le p
\le \infty$ and arguing as above, we know that the map $\Lambda:
L_p(\Om, L_p(\mathcal{M})) \to \ell_p(L_p(\mathcal{M}))$ defined
by $\Lambda(f) = (\int f \eps_k d\mu)_{k \le n}$ is a complete
contraction. This yields the result for $p \ge 2$. When $p<2$ we
note that $\Lambda: L_\infty(\Omega;L_1(\mathcal{M})) \to
\ell_{\infty}(L_1(\mathcal{M}))$ is completely bounded. Again
interpolation yields the result. The sharpness of this result is
justified in Remark \ref{Sharpness} below. \fin

\begin{lemma} \label{cintpol} Let $A$ be a
$C^*$-algebra and $\phi$ be a state on $\mathcal{N} = A^{**}$
whose restriction to $A$ is faithful. Let $d \in L_1(\mathcal{N})$
be the associated density with support $e$ in $\mathcal{N}$. Let
us set $\mathcal{N}_e = e \mathcal{N} e$. Then, we have
$$\big[ A, L_1(\mathcal{N}_e) \big]_{\frac1p} =
L_p(\mathcal{N}_e).$$
\end{lemma}

\dem Following Kosaki's work, we have symmetric injective
embedding of $\mathcal{N}_e$ in $L_1(\mathcal{N}_e)\cong
(\mathcal{N}_e)^*$ given by $\iota(x)=d^{1/2}xd^{1/2}$. Let $x\in
A$ and $y$ be an analytic element in $\mathcal{N}_e$. Then we
observe that
 \[ \langle \iota(x),y\rangle = tr(d^{1/2}xd^{1/2}y)
 =tr(dxd^{1/2}yd^{-1/2})=\phi(x\sigma_{i/2}(y)) \, .\]
Since the elements of the form $\sigma_{i/2}(y)$, $y$ analytic,
are in dense in $\mathcal{N}_e$, we deduce from $\iota(x)=0$ that
$\phi(x^*x)=0$. However, $\phi$ is faithful and hence
$(A,L_1(\mathcal{N}_e))$ is indeed an interpolation couple and we
may define $X_p = [A, L_1(\mathcal{N}_e)]_{1/p}$. By Kaplansky's
density theorem we known that $eAe$ is strongly dense in
$\mathcal{N}_e$ and hence $d^{1/2}Ad^{1/2}$ norm dense in
$L_1(\mathcal{N}_e)$. Thus the interpolation couple has dense
intersection. The unit ball in $X_p^*$ is the closure in the sum
topology of the unit ball in
$$Z_p = \big[ A^*, L_1(\mathcal{N}_e)^* \big]_{\frac1p} = \big[
\mathcal{N}_e, L_1(\mathcal{N}_e) \big]_{1-\frac1p},$$ see
\cite{BL} for further details. Here the natural inclusion map is
again given by
 $$n \in \mathcal{N}_e \mapsto d^\frac12 n d^\frac12 \in
 L_1(\mathcal{N}_e) \mapsto d^\frac12 n d^\frac12|_{A} \in A^*,$$
because $A$ is the  intersection in the interpolation couple.
Certainly, $L_1(\mathcal{N}_e)$ is faithfully embedded in $A^*$.
Thus in the dual picture we find exactly the symmetric version of
Kosaki's embedding  \cite{Ko}, $Z_p = L_{p'}(\mathcal{N}_e)$.
Since $L_{p'}(\mathcal{N}_e)$ is reflexive, its unit ball is
already closed in the sum topology. Indeed, given any converging
sequence in the sum topology, it is easily checked that the limit
is a cluster point of the sequence in the weak$^*$ topology. This
gives $X_p^*= L_{p'}(\mathcal{N}_e)$, so that the inclusion $X_p
\subset L_p(\mathcal{N}_e)$ is isometric. The assertion then
follows from the fact that the norm dense subspace $d^{1/2p} A \,
d^{1/2p}$ of $L_p(\mathcal{N}_e)$ is contained in $X_p$. \fin

\begin{lemma} \label{uf} Let $\U$ be an ultrafilter
on an index set $I$ and $$(d_i)^\bullet \in \prodd_\U
L_1(\mathcal{M}).$$ Let $\phi(x) = \lim_{i,\U} \mathrm{tr}(d_ix)$
be the corresponding weak limit state and $d \in
L_1(\mathcal{M}^{**})$ the corresponding nonfaithful density
supported by $e$ in $\mathcal{M}^{**}$. Then, there exists a
completely contractive map densely defined on $d^{1/2p}
\mathcal{M} \, d^{1/2p}$ by
$$u_p: d^{\frac{1}{2p}} x \, d^{\frac{1}{2p}} \in
L_p(e\mathcal{M}^{**}e) \mapsto \big( d_i^{\frac{1}{2p}} x \,
d_i^{\frac{1}{2p}} \big)^\bullet \in \prodd_\U L_p(\mathcal{M}).$$
\end{lemma}

\dem Let $e_\U$ be the support of $$\phi_\U (x_i)^\bullet =
\limm_{i,\U} \mathrm{tr}(d_ix_i)$$ and consider the
$\sigma$-finite von Neumann algebra $\mathcal{M}_\U = e_\U \big[
\prod_\U L_1(\mathcal{M}) \big]^* e_\U$. The image of $u_p$ sits
on $L_p(\mathcal{M}_\mathcal{U})$ and $\phi_\mathcal{U}$ is
faithful on $\mathcal{M_U}$. Hence, the spaces
$L_p(\mathcal{M}_\U)$ interpolate by Kosaki's result. Let $f$ be
the support of $\phi$ in ${\mathcal M}$ and $e$ be the support of
$\phi$ in ${\mathcal M}^{**}$. Note that $e\le f$. We apply Lemma
\ref{cintpol} to $A= \mathcal{M}_f = f{\mathcal M}f$ and obtain
$$L_p(e\mathcal{M}^{**}e) = \big[ \mathcal{M}_f,
L_1(e \mathcal{M}^{**} e)\big]_{\frac1p}.$$ Therefore, the map
$u_p$ is obtained by interpolation. Clearly $u_\infty(x)= e_\U
(x)^\bullet e_\U$ is a complete contraction. The interesting part
is the case $p=1$. For a positive $x \in \mathcal{M}$, we note
that
$$\big\| (d_i^{\frac12} x \, d_i^{\frac12} )^\bullet
\big\|_{\prod_\U L_1(\mathcal{M})} = \limm_{i,\U} \mathrm{tr}
\big( d_i^{\frac12} x \, d_i^{\frac12} \big) = \limm_{i,\U}
\mathrm{tr}(d_ix) = \phi(x).$$ For a positive element $x \in
\mathcal{M}^{**}$, we may apply Kaplansky's density theorem and
approximate $x^{1/2}$ in $\mathrm{SOT} \cap \mathrm{SOT}^*$ by a
net $x_{\la} \in \mathcal{M}$ such that $\|x_{\la}\| \le
\|x^{1/2}\|$. Then we have
 $$\limm_{i,\U} \big\| (x_{\la} - x_{\mu}) d_i^{\frac12} \big\|_2^2
 = \limm_{i,\U} \mathrm{tr} \big( d_i |x_{\la} - x_{\mu}|^2 \big) =
 \phi \big( |x_{\la}-x_{\mu}|^2 \big).$$
Hence, $(x_{\la} d_i^{\frac12})^{\bullet}$ is Cauchy in $\prod_\U
L_2(\mathcal{M})$ with limit $(x^{\frac12}
d_i^{\frac12})^{\bullet}$ because
\begin{eqnarray*}
(\phi \big( |x_{\la}-x_{\mu}|^2 \big))^{1/2} & \le & \phi \big(
|x_{\la}-x|^2 \big)^\frac12 + \phi \big( |x_{\mu}-x|^2 \big)^\frac12 \\
& = & \big( \phi(x_{\la}^*x_{\la})+\phi(x^*x)-\phi(x_{\la}^*x)-
\phi(x^*x_{\la}) \big)^\frac12 \\ & + & \big(
\phi(x_{\mu}^*x_{\mu})+\phi(x^*x)-\phi(x_{\mu}^*x)-
\phi(x^*x_{\mu}) \big)^\frac12
\end{eqnarray*}
converges to 0. Moreover, we have
 \[ u_1(x) = (d_i x^{\frac12})^{\bullet} (x^{\frac12} d_i)^{\bullet} \in
 \prodd_\U L_1(\mathcal{M}).\]
Now let $x \in e \mathcal{M}^{**} e$ be a self-adjoint element.
Then we recall from \cite{HJX} that
$$\limm_{i,\U} \big\| d_i^{\frac12} x \, d_i^{\frac12} \big\|_1
\le \inf_{x=x_1-x_2} \phi(x_1) + \phi(x_2) = \big\| d^{\frac12} x
\, d^{\frac12} \big\|_1$$ where the infimum is taken over positive
elements in $e \mathcal{M}^{**} e$. This implies that $$u_1:
d^{\frac12} x \, d^{\frac12} \in L_1(e \mathcal{M}^{**} e) \to
\big( d_i^{\frac12} x \, d_i^{\frac12} \big)^\bullet \in \prodd_\U
L_1(\mathcal{M}),$$ is a c.p. map with $u_1^*(\mathbf{1})=
\mathbf{1}$. Hence, $u_1$ and $u_1^*$ are contractions.
Interpolation and the density of
$[\mathcal{M}_f,L_1(e\mathcal{M}^{**} e)]_{1/p} \subset
L_p(e\mathcal{M}^{**}e)$ implies the result. Since the same
argument holds for $M_m(e\mathcal{M}^{**} e)$, $u_p$ is a complete
contraction. \fin

\demA Let us begin by proving the statement i). Let $\mathcal{N} =
A^{**}$ and consider the adjoint mapping $T^*: X^* \to A^*$. Since
$A^* \simeq_{cb} L_1(\mathcal{N}^{\mathrm{op}})$ and $T$ is a
completely $(p,1)$-summing map we deduce from (the dual version
of) Remark \ref{linn} that $\|T^* \ten id: \ell_{p'}(X^*) \to
L_1(\mathcal{N}^{\mathrm{op}}; \ell_{\infty})\|_{cb} \le
\pi_{p,1}^{cb}(T)$. According to Theorem B, this implies
$$\big\| T^* \ten id: S_{q'}(X^*) \to L_1(\mathcal{N}^{\mathrm{op}};
S_{q'}) \big\|_{cb} \le  c(p,q) \, \pi_{p,1}^{cb}(T).$$ Dualizing
again, we obtain the following key inequality $$\big\| T \ten id:
A(S_q) \to S_q(X) \big\|_{cb} \le c(p,q) \, \pi_{p,1}^{cb}(T).$$
Here we interpret $A(S_q)$ as in Remark \ref{linn}
\begin{eqnarray*}
A(S_q^m) & = & \big[ M_m(A), A(S_1^m) \big]_{\frac1q}
\\ & \subset & \big[ L_\infty(A^{**}; S_\infty^m), L_\infty(A^{**}; S_1^m)
\big]_{\frac1q} \ = \ L_\infty(\mathcal{N}; S_q^m).
\end{eqnarray*}
Now we follow Pisier and apply the Grothendieck-Pietsch separation
argument as in \cite[Theorem 5.1]{P2}. Namely, the substitute of
the auxiliary Theorem 5.3 there for $A \ten_\mathrm{min} S_q$ has
to be replaced here by the fact that $A(S_q) \subset
L_\infty(\mathcal{N}; S_q)$ is understood as a conditional
$L_\infty$ space with norm given by $$\|x\|_{A(S_q)} = \sup \Big\{
\|\alpha x \beta \|_{L_q(\mathcal{N} \bar\ten
\mathcal{B}(\ell_2))} \, \big| \ \|\alpha\|_{L_{2q}(\mathcal{N})},
\|\beta\|_{L_{2q}(\mathcal{N})} \le 1\Big\}.$$ Then, it turns out
that Pisier's argument in \cite{P2} generalizes verbatim to this
setting and we find nets $(a_{\la})$ and $(b_{\la})$ in the
positive part of the unit ball of $L_{2q}(\mathcal{N})$ satisfying
the inequality $$\|T(x)\|_{S_q(X)}\le c(p,q) \, \pi_{p,1}^{cb}(T)
\, \limm_{\la} \big\| a_{\la} x b_{\la}
\big\|_{S_q(L_q(\mathcal{N}))}.$$ On $M_2(\mathcal{N})$, we define
the state $$\phi(x) = \limm_{\la} \frac12 \Big[ \mathrm{tr} \big(
a_{\la}^{2q} x_{11} \big) + \mathrm{tr} \big( b_{\la}^{2q} x_{22}
\big) \Big].$$ Let $d \in L_1(M_2(\mathcal{N}))$ be the density of
$\phi$. We also use the notation $d_a, d_b$ for the densities of
the states $\phi_a(x) = \lim_{\la} \mathrm{tr}(a_{\la}^{2q} x)$
and $\phi_b(x) = \lim_{\la} \mathrm{tr}(b_{\la}^{2q} x)$.
According to Lemma \ref{uf}, we see that $$u_q \big(
d^{\frac{1}{2q}} x \, d^{\frac{1}{2q}} \big) = \big(
d_{\la}^{\frac{1}{2q}} x \, d_{\la}^{\frac{1}{2q}} \big)^{\bullet}
\quad \mbox{with} \quad d_\lambda = \frac12 \Big( e_{11} \ten
a_{\la}^{2q} + e_{22} \ten b_{\la}^{2q} \Big).$$ is a complete
contraction. Restricting this to the $(1,2)$ entry, we deduce that
\begin{eqnarray*}
\limm_{\la} 2^{-\frac1q} \big\| a_{\la} x b_{\la} \big\|_q & = &
\limm_{\la} \big\| d_\la^{\frac{1}{2q}} (e_{12} \ten x) \,
d_\la^{\frac{1}{2q}} \big\|_q \\ & \le & \big\| d^{\frac{1}{2q}}
(e_{12} \ten x) \, d^{\frac{1}{2q}} \big\|_q \ = \ 2^{-\frac1q}
\big\| d_a^{\frac{1}{2q}} x \, d_b^{\frac{1}{2q}} \big\|_q.
\end{eqnarray*}
Moreover, the same chain of inequalities holds for $x$ replaced by
an element in $M_m(A)$. The first assertion then follows
immediately. Indeed, it just remains to choose the densities
$\delta_1^{2q} = d_a$, $\delta_2^{2q} = d_b$ and define the map
$$w(\delta_1 x \delta_2) = T(x).$$ To prove ii), we follow an
argument by Haagerup. According to the first part applied to $A =
\mathcal{M}$, we find $\delta_1, \delta_2 \in
L_{2q}^+(\mathcal{M}^{**})$ of norm $1$. Moreover, using the
existence of a central projection $z$ in $\mathcal{M}^{**}$ such
that $\mathcal{M} = z \mathcal{M}^{**}$, we define $d_1 = z
\delta_1$ and $d_2 = z \delta_2$. Let $(z_{\la}) \subset
\mathcal{M}$ be a net of contractions which converges strongly to
$z$ in $\mathcal{M}^{**}$. Then $d_1 = \lim_{\la} z_{\la}
\delta_1$ and $d_2 = \lim_{\la} z_{\la} \delta_2$. On the other
hand $\mathbf{1}_\mathcal{M}-z_{\la}$ converges strongly to $0$,
where strongly refers this time to $\mathcal{M}$. Since $T$ is
supposed to be normal, we have $T^*(X^*) \subset
L_1(\mathcal{M})$. This implies $$\limm_{\la,\mu} \big\langle x^*,
T(z_{\la}xz_{\mu}) \big\rangle = \limm_{\la,\mu} \big\langle
z_{\la} T^*(x^*) z_{\mu} , x \big\rangle =  \langle x^*, T(x)
\rangle.$$ Let $x \in S_q(\mathcal{M})$ and $x^*$ in the unit ball
of $S_{q'}(X^*)$ so that
$$\|T(x)\|_{S_q(X)} = \big| \langle x^*, T(x) \rangle \big|.$$
Then we find
\begin{eqnarray*}
\|T(x)\|_{S_q(X)} & = &  \limm_{\la,\mu} \big| \big\langle x^*,
T(z_{\la}xz_{\mu}) \big\rangle \big| \\
& \le & c(p,q) \, \pi_{p,1}^{cb}(T) \limm_{\la,\mu} \big\|
\delta_1 z_{\la} x z_{\mu} \delta_2
\big\|_{S_q(L_q(\mathcal{M}^{**}))}
\\ & = & c(p,q) \, \pi_{p,1}^{cb}(T) \big\| d_1 x d_2
\big\|_{S_q(L_q(\mathcal{M}^{**}))}.
\end{eqnarray*}
This shows that $v(d_1 x d_2) = T(x)$ is continuous and even
completely bounded. The proof of iii) follows the same pattern
above. We first dualize and consider the map $T^*: X^* \to
L_{s'}(\mathcal{M})$, which is $(p',\infty)$-convex. Indeed, this
follows by duality since
$$L_s(\mathcal{M}; \ell_1) \stackrel{id}{\longrightarrow} \ell_1
\ten_{\mathrm{min}} L_s(\mathcal{M}) \stackrel{T}{\longrightarrow}
\ell_p(X)$$ is completely bounded. Then, since $s' < q' < p'
\wedge \infty$, we may apply Theorem B to deduce that $T^* \ten
id: S_{q'}(X^*) \to L_{s'}(\mathcal{M}; S_{q'})$ is completely
bounded with cb-norm controlled by $c(p,q,s) \,
\pi_{p,1}^{cb}(T)$. Dualizing back and with the help of the
Grothendieck-Pietsch factorization theorem (adapted to this
setting as indicated above), we find nets $(a_{\la})$, $(b_{\la})$
in the positive part of the unit ball of $L_{2w}(\mathcal{M})$
such that
$$\|T(x)\|_{S_q(X)} \le c(p,q,s) \, \pi_{p,1}^{cb}(T) \,
\limm_{\la} \big\| a_{\la} x b_{\la}
\big\|_{S_q(L_q(\mathcal{M}))}.$$ Let us assume for simplicity
that $a_{\la} = b_{\la} = d_\la$. Recall that this can always be
done using the $2\times 2$ matrix trick from above. Then we define
the following weak$^*$ limit in $L_{s'}(\mathcal{M})$
$$\mathrm{tr}(dx)  = \limm_{\la} \mathrm{tr} \big( d_{\la}^{2w/s'}
x \big).$$ The assertion is obtained from the inequality
$$\limm_{\la} \big\| d_{\la} x \hskip1pt d_{\la} \big\|_q \le
\big\| d^{s'/2w} x \, d^{s'/2w} \big\|_q,$$ which follows by
approximating $x \sim d_\la^{w/s} z d_\la^{w/s}$ and applying
Lemma \ref{uf}. \fin

\begin{remark}
\emph{According to Remark \ref{Constante}, we obtain
$$c(p,q) \lesssim \frac{1}{1 -
\frac{p}{q}}.$$ We also have the weaker estimate $c(p,q,s)
\lesssim q(s-p)/(q-p)$ for $s < \infty$.}
\end{remark}

\begin{remark} \label{Splpsumming}
\emph{We may define canonically $$\pi_{p,q}^{cb}(T) = \big\| id
\ten T: \ell_q \ten_{\min} X \to \ell_p(Y) \big\|_{cb}$$ as the
completely $(p,q)$-summing norm of $T: X \to Y$. At the time of
this writing, it is not clear whether $\pi_{p,p}^{cb}(T) =
\pi_p^o(T)$ holds for all maps $T$. However, Pisier's
factorization theorem immediately implies that every completely
$p$-summing map is completely $(p,p)$ summing, and
$\pi_{p,p}^{cb}(T) \le \pi_p^o(T)$. If in addition $T$ is a normal
map on an injective von Neumann algebra, then the norms are
equivalent. Indeed, let $T^*: Y^* \to L_1(\mathcal{M})$ be the
adjoint, $\mathcal{M}$ injective such that
$$\pi_{p,p}^{cb}(T) = \big\| id \ten T^*: \ell_{p'}(Y^*) \to
L_1(\mathcal{M}; \ell_{p'}) \big\|_{cb} < \infty.$$ Recall from
\cite{JP1} that we have a cb-embedding $j: S_{p'}^m \to
L_1(\mathcal{N}; \ell_{p'})$, so that $$id \ten j:
L_1(\mathcal{M}; S_{p'}^m) \to L_1(\mathcal{M} \bar\ten
\mathcal{N}; \ell_{p'})$$ is an isomorphic embedding. This map
uses independent copies and hence it is easy to check that $j\ten
id_{Y^*}:S^m_{p'}(Y^*)\to L_1(\mathcal{N}; \ell_{p'}(Y^*))$
remains bounded with a constant $c(p)$. Then we find the following
diagram
\begin{align*}
\begin{array}{ccc}
S_{p'}^m(Y^*) &  \stackrel{T^*}{\longrightarrow} &
L_1(\mathcal{M},S_{p'}^m) \\ j \downarrow & & \uparrow j^{-1} \\
L_1(\mathcal{N}; \ell_{p'}(Y^*)) & \stackrel{T^*}{\longrightarrow}
& L_1(\mathcal{M} \bar\ten \mathcal{N}; \ell_{p'}).
\end{array}
\end{align*}
The two maps $\downarrow$ and $\uparrow$ are bounded, and hence
$$\big\| id_{S_{p'}^m} \ten T^*: S_{p'}^m(Y^*) \to
L_1(\mathcal{M}; S_{p'}^m) \big\| \le c(p)\, \pi_{p,p}^{cb}(u)$$
is still bounded with constants independent of $m$. This completes
the argument.}
\end{remark}

\subsection{Applications I. Operator spaces} Our first application
is an operator space analog of Rosenthal's theorem \cite{Ro} for
subspaces of (commutative or not) $L_p$ spaces. This partly
justifies our definition of cb-cotype, see \cite{GP2,JP1,Lee,Pa}
for related notions.

\demCA We shall prove
i)$\Rightarrow$ii)$\Rightarrow$iii)$\Rightarrow$i). The first
implication follows from Lemma \ref{prel}. For the second
implication, assume that $X^*$ is completely $(p_0',1)$-summing
for some index $p < p_0 < 2$ and let $j: X \to L_p(\mathcal{M})$
be the inclusion map. Take the (necessarily normal) adjoint map $T
= j^*: L_{p'}(\mathcal{M}) \to X^*$. Given $p_0' < q' < p'$, the
map $T: \ell_1 \ten_\mathrm{min} L_{p'}(\mathcal{M}) \to
\ell_{p_0'}(X^*)$ is completely bounded since $id_{X^*}$ is
completely $(p_0',1)$-summing and $$\ell_1 \ten_\mathrm{min}
L_{p'}(\mathcal{M}) \stackrel{T}{\longrightarrow} \ell_1
\ten_\mathrm{min} X^* \stackrel{id}{\longrightarrow}
\ell_{p_0'}(X^*).$$ In particular, $T$ satisfies the assertion of
Theorem A. Let $v: L_{q'}(\mathcal{M}) \to X^*$ be the
corresponding map. Then $v^*: X \to L_q(\mathcal{M})$ is also
completely bounded and $d_1 v^*(x) \, d_2 = j(x)$. In particular,
since $d_1, d_2$ are norm 1 in $L_{2w}(\mathcal{M})$ and $\frac1p
= \frac1q + \frac1w$
$$\|x\|_{M_m(X)} = \|j(x)\|_{M_m(L_p(\mathcal{M}))} = \big\| d_1
v^*(x) \, d_2 \big\|_{M_m(L_p(\mathcal{M}))} \le
\|v^*(x)\|_{M_m(L_q(\mathcal{M}))}.$$ Thus, $X$ is cb-isomorphic
to $v^*(X) \subset L_q(\mathcal{M})$. For the third implication,
the Rademacher transform map $\Lambda: f \in
\Rad(L_{q'}(\mathcal{M})) \mapsto (\int_\Omega f \, \eps_k \,
d\mu) \in \ell_{q'}(L_{q'}(\mathcal{M}))$ is completely
contractive and this remains true for every quotient of
$L_{q'}(\mathcal{M})$. In particular, $X^*$ has cb-cotype $q'$.
The proof is complete. \fin

\begin{corollary} \label{gen}
If $p \ge 2$ and $id_X$ is completely $(p,1)$-summing
$$\Pi_{1}^o(X,Y)=  \Pi_{q'}^o(X,Y) \quad \mbox{for all operator
spaces $Y$ and $q>p$}.$$
\end{corollary}

\dem The inclusion $$\Pi_1^o(X,Y) \subset \Pi_{q'}^o(X,Y)$$ is
well-known. For the converse, we consider $u: M_m \to X$ and note
that $$\pi_{p,1}^{cb}(u)\le \|u\|_{cb} \, \pi_{p,1}^{cb}(id_X).$$
Theorem \ref{tha} for $\mathcal{M}=M_m$ gives $a,b \in S_{2q}^m$
and a cb-map $w: S_q^m \to X$ such that $$u = w \circ M_{ab} \quad
\mbox{and} \quad \|a\|_{2q} \|w\|_{cb} \|b\|_{2q} \le c(p,q) \,
\|u\|_{cb} \, \pi_{p,1}^{cb}(id_X).$$ The argument now follows by
a standard duality argument. We refer the reader to \cite[Chapter
7]{P2} for a brief review of the duality theory of $p$-summing
maps both in the Banach and operator space settings. We shall also
use the $p$-nuclear norm $\nu_p^o$ and the fact that it is trace
dual to $\pi_{q'}^o$, see \cite[Chapter 3]{J0}. If $T: X \to Y$
and $v: Y \to M_m$, we deduce that
\begin{eqnarray*}
\big| \mathrm{tr}(vTu) \big| & = & \big| \mathrm{tr}(M_{ab}vTw)
\big| \\ & \le & \nu_{q}^o(M_{ab}v) \pi_{q'}^o(Tw) \\
& \le & \|v\|_{cb} \, \|a\|_{2q} \|w\|_{cb} \|b\|_{2q} \,
\pi_{q'}^o(T) \\ & \le & c(p,q) \pi_{p,1}^{cb}(id_X) \, \|v\|_{cb}
\|u\|_{cb} \, \pi_{q'}^o(T).
\end{eqnarray*}
Thus we obtain the inequality
$$\sup_{\|u\|_{cb}, \|v\|_{cb} \le 1} \big| \mathrm{tr}(vTu) \big|
\le c(p,q) \, \pi_{p,1}^{cb}(id_X) \, \pi_{q'}^o(T).$$ Since
$\mathcal{CB}(Y,M_m) = [S_1^m(Y)]^*$ and $\mathcal{CB}(M_m,X) =
S_1^m \ten_{\min} X$, we get $$\big\| T \ten id: S_1^m \ten_{\min}
X \to S_1^m(Y) \big\| \le c(p,q) \, \pi_{p,1}^{cb}(id_X) \,
\pi_{q'}^o(T),$$ but the left hand side is the completely
$1$-summing of $T$. The proof is complete. \fin

\demCB The first assertion follows from Theorem A, while the
second assertion follows from Lemma \ref{elm} applied to $OH$ and
Corollary \ref{gen}. \fin

\subsection{Applications II. Noncommutative $L_p$ spaces}

We now investigate some further consequences of our results for
linear maps between noncommutative $L_p$ spaces equipped with
their natural operator space structures.

\begin{corollary}
Let $2 \le q_1 < p_1 < q_2 \le p_2 \le \infty$. Assume that $$T:
L_{p_2}(\mathcal{M}) \to L_{p_1}(\mathcal{M}) \quad \mbox{and}
\quad S: L_{q_2}(\mathcal{N}) \to L_{q_1}(\mathcal{N})$$ are
completely bounded maps with $\mathcal{M}, \mathcal{N}$ being
$\mathrm{QWEP}$ von Neumann algebras. In the case $p_2=\infty$ or
$q_2=\infty$, assume in addition that the corresponding map is
normal. Then, the following map is completely bounded $$T \ten S:
L_{p_2} \big( \mathcal{M}; L_{q_2}(\mathcal{N}) \big) \to L_{p_1}
\big( \mathcal{M}; L_{q_1}(\mathcal{N}) \big).$$
\end{corollary}

\dem If $2 \le p_1 < q_2 \le p_2$, we claim that
$$\big( \widetilde{T \ten id} \big) (x\ten y)=y\ten T(x)$$ satisfies
$$\big\| \widetilde{T \ten id}: L_{p_2} \big( \mathcal{M};
L_{q_2}(\mathcal{N}) \big) \to L_{q_2} \big( \mathcal{N};
L_{p_1}(\mathcal{M}) \big) \big\|_{cb} \le c(p_1,q_2) \,
\|T\|_{cb}.$$ Indeed, if $p_2=q_2$ then $L_{p_2}(\mathcal{M};
L_{p_2}(\mathcal{N})) = L_{p_2}(\mathcal{M} \bar\ten \mathcal{N})
= L_{p_2}(\mathcal{N}; L_{p_2}(\mathcal{M}))$. Since $\mathcal{N}$
is QWEP, we deduce the assertion from the complete boundedness of
$T$. A similar argument can be found in \cite{J4}. When $p_2 >
q_2$, we use that $L_{p_1}(\mathcal{M})$ has cb-cotype $p_1$ and
Theorem A to factorize $$T = v \circ M_{ab},$$ where $v:
L_{q_2}(\mathcal{M}) \to L_{p_1}(\mathcal{M})$ is completely
bounded and $M_{ab}(x)=axb$ with $a,b$ positive norm 1 elements of
$L_{2s}(\mathcal{M})$ for $1/q_2=1/p_2+1/s$. It is clear that the
map $$M_{ab} \ten id: L_{p_2}(\mathcal{M}; L_{q_2}(\mathcal{N}))
\to L_{q_2}(\mathcal{M}; L_{q_2}(\mathcal{N}))$$ is completely
contractive. Moreover, our argument for $p_2=q_2$ gives
\begin{eqnarray*}
\big\| v \ten id: L_{q_2} \big( \mathcal{M}; L_{q_2}(\mathcal{N})
\big) \to L_{q_2} \big( \mathcal{N}; L_{p_1}(\mathcal{M}) \big)
\big\|_{cb} & \le & c(p_1,p_2,q_2) \, \pi_{p_1,1}^{cb}(T) \\
& \le & c(p_1,p_2,q_2) \, \|T\|_{cb}.
\end{eqnarray*}
This proves our claim. Moreover, if $$2 \le q_1 < p_1 < q_2 \le
p_2$$ the same argument for $\widetilde{id\ten S}$ yields
$$\big\| \widetilde{id \ten S}: L_{q_2} \big( \mathcal{N};
L_{p_1}(\mathcal{M}) \big) \to L_{p_1} \big( \mathcal{M};
L_{q_1}(\mathcal{N}) \big) \big\|_{cb} \le c(p_1,q_1,q_2) \,
\|S\|_{cb}.$$ Combining the two estimates, we deduce the
assertion. The proof is complete. \fin

\begin{corollary} \label{1p}
If $2 \le p < q < \infty$ and $\mathcal{M}, \mathcal{N}$ are
hyperfinite $$\mathcal{CB} \big( L_1(\mathcal{M}),L_p(\mathcal{N})
\big) = \Pi_q^o \big( L_1(\mathcal{M}),L_p(\mathcal{N}) \big).$$
\end{corollary}

\dem Since $1 < p' \le 2$ and according to \cite{JP3,JP4}, we have
a cb-embedding $j: L_{p'}(\mathcal{N}) \to L_1(\mathcal{A})$ for
some hyperfinite von Neumann algebra $\mathcal{A}$. The dual map
$j^*: \mathcal{A}^{op} \to L_p(\mathcal{N})$ is a complete
surjection. Let $u: L_1(\mathcal{M}) \to L_p(\mathcal{N})$ be a
completely bounded map and $u^*: L_{p'}(\mathcal{N})\to
\mathcal{M}^{op}$ its adjoint map. Since $\mathcal{M}$ is
injective, we have a cb-norm preserving extension
$w:L_1(\mathcal{A})\to \mathcal{M}^{op}$. The restriction
$\tilde{u}$ of $w^*(\mathcal{M}^{op})^*\to \mathcal{A}^{op}$ to
$L_1(\mathcal{M})$ gives an extension of $u:L_1(\mathcal{M})\to
\mathcal{A}^{op}$ such that $u=j^*w^*$ and
$$\|\tilde{u}\|_{cb} \le \|u\|_{cb} \|j\|_{cb} \|j^{-1}\|_{cb} \le c
\|u\|_{cb}.$$ Since $L_p(\mathcal{M})$ has cb-cotype $p$ from
Lemma \ref{elm} and $j^*$ is normal, we know from Theorem A that
$j^*$ is completely $q$-summing. Recall that the fact that
$\mathcal{A}$ is injective is used here to ensure that
$\mathcal{A}(S_q^m) = \mathcal{A} \ten_{\min} S_q^m$. Thus we
conclude $u = j^* \tilde{u}$ is also completely $q$-summing. \fin

\begin{corollary}
If $\mathcal{M}$ is finite and hyperfinite and $$T:
L_1(\mathcal{M}) \to L_2(\mathcal{M})$$ is completely bounded,
then the eigenvalues of $T: L_2(\mathcal{M}) \to L_2(\mathcal{M})$
satisfy $$\Big( \summ_k |\la_k(T)|^2 \Big)^\frac12 \le
\|T\|_{cb}.$$
\end{corollary}

\dem It is well-known \cite[3.4.3.13]{J0} that $$\Big( \summ_k
|\la_k(T)|^q \Big)^\frac1q \le \pi_q^o(T)$$ for $2 < q < \infty$.
Here $\la_k(T)$ are the eigenvalues in non-decreasing order. Let
us take the opportunity to correct an oversight in the proof. In
\cite[p.238]{J0} it was claimed that \[ \prodd_\U S_p
\stackrel{?}{=} \Big[ \prodd_\U S_{\infty}, \prodd_\U S_2
\Big]_{\frac2p}\] interpolates. However, this is not an
interpolation couple. Instead, one has to use Pisier's
factorization theorem and use that for a positive density
$d=(d_i)\in \prod_\U S_1$ the spaces $$LS_p = \mathrm{cl} \Big(
\Big\{ (d_i^{1/2p}x_id_i^{1/2p})^{\bullet} \, \big| \,
(x_i)^{\bullet}= e_\U(x_i)^{\bullet}e_\U \in (\prodd_{\U} S_1)^*
\Big\} \Big) \subset \prodd_\U S_p$$ form an interpolation scale,
due to Kosaki's interpolation theorem. In the rest of the proof
one works with these spaces. In order to push the result to $q=2$,
we may apply a standard tensor trick. Let $m \in \N$ and $j_m:
L_2(\mathcal{M}^{\ten_m}) \to L_1(\mathcal{M}^{\ten_m})$ be the
natural completely contractive inclusion map. Then we deduce from
Corollary \ref{1p} that
\begin{eqnarray*}
\Big( \sum_{k=1}^n |\la_k(T)|^2 \Big)^{\frac{m}{2}}
& = & \Big( \sum_{k=1}^{n^m} |\la_k(T^{\ten_m})|^2 \Big)^\frac12 \\
& \le & n^{\frac{m}{2}-\frac{m}{q}} \Big( \sum_{k=1}^{n^m}
|\la_k(T^{\ten_m})|^q \Big)^{\frac1q} \\
[7pt] & \le & n^{\frac{m}{2}-\frac{m}{q}} \, \pi_q^o(T^{\ten_m}
j_m) \ \le \ n^{\frac{m}{2}-\frac{m}{q}} c(q) \, \|T^{\otimes
m}\|_{cb}.
\end{eqnarray*}
We now claim that $\|T^{\otimes m}\|_{cb} \le \|T\|_{cb}^m$.
Indeed, given $$T: L_1(\mathcal{M}) \to L_2(\mathcal{M}) \quad
\mbox{and} \quad S: L_1(\mathcal{N}) \to L_2(\mathcal{N}),$$ we
observe that $$L_1(\mathcal{M} \bar\otimes \mathcal{N})
\stackrel{S}{\longrightarrow} L_1(\mathcal{M}; L_2(\mathcal{N}))
\longrightarrow L_2(\mathcal{N}; L_1(\mathcal{M}))
\stackrel{T}{\longrightarrow} L_2(\mathcal{M} \bar\otimes
\mathcal{N})$$ where the middle map is a complete contraction by
Minkowski's inequality. Hence we have $\|T \otimes S\|_{cb} \le
\|T\|_{cb} \|S\|_{cb}$. Applying it $m-1$ times, we deduce our
claim and therefore we get $$\Big( \sum_{k=1}^n |\la_k(T)|^2
\Big)^{\frac{m}{2}} \le n^{\frac{m}{2}-\frac{m}{q}} c(q) \,
\|T\|_{cb}^m.$$ Thus, taking $m$-th roots and sending $(m,q) \to
(\infty,2)$, the result follows. \fin

\begin{corollary} \label{eig}
Let $1 < p < \infty$ and $q > p \vee p'$. If $\mathcal{M}$ is
hyperfinite and the map $T: \mathcal{M} \to \mathcal{M}$ is normal
with a factorization $T=vw$, where $v: L_p(\mathcal{M}) \to
\mathcal{M}$ and $w: \mathcal{M} \to L_p(\mathcal{M})$ normal,
both completely bounded. Then, we have $$\Big( \summ_k
|\la_k(T)|^q \Big)^\frac1q \le c(p,q) \, \|v\|_{cb} \|w\|_{cb}.$$
\end{corollary}

\dem When $p \ge 2$, this follows from Theorem A because $w:
\mathcal{M} \to L_p(\mathcal{M})$ is completely $q$-summing and
hence $T=vw$ is also completely $q$-summing. In the case $1<p<2$,
we consider $T^*=w^* v^*$ and deduce from Corollary \ref{1p} that
the map $v^*: L_1(\mathcal{M}) \to L_{p'}(\mathcal{M})$ is
completely $q$-summing. Following the eigenvalue estimates from
\cite[p.238]{J0} and letting $T_* = T^*|_{L_1}$, we know that
$T_*^m$ is compact for some $m\in \mathbb{N}$. Hence $T^m$ is also
compact and $T$ is a Riesz operator. Recall that an operator
$T:X\to X$ is Riesz if for all $\eps>0$ there exist $n,m\in
\mathbb{N}$ and $y_1,...,y_m$ such that $T^n(B_{X})\subset
\bigcup_k y_k+\eps B_{X}$. Fortunately, we know by a result of
West which can be found in \cite[3.2.26]{Pi} that for a Riesz
operator the eigenvalues sequence  $(\la_k(T))$ can be arranged so
that $(\la_k(T))=(\la_k(T^*))$. We also refer to \cite[Section
3.2]{Pi} for the definition of the eigenvalue sequence respecting
the algebraic multiplicity. Hence our estimate of the $q$-norm of
$(\la_k(T_*))$ implies the same estimate for the eigenvalue
sequence of $T:\mathcal{M}\to \mathcal{M}$. \fin

\begin{remark} \label{Sharpness}
\emph{Let us consider an example. Given a sequence $(\mu_k) \in
\ell_p$ of positive numbers, the cb-norm of the diagonal map
$\Delta_{\sqrt{\mu}}: e_{k1} \in C \mapsto \sqrt{\mu_k} \, e_{k1}
\in C_p$ is given by
$$\big\| \Delta_{\sqrt{\mu}}: C \to C_p \big\|_{cb} = \Big(
\sum_{k=1}^\infty \mu_k^p \Big)^{\frac{1}{2p}} = \big\|
\Delta_{\sqrt{\mu}}: C_p \to C \big\|_{cb}.$$ Hence,
$\Delta_{\sqrt{\mu}}$ factors through $S_p$ and $S_{p'}$ and
therefore the best possible exponent in Corollary \ref{eig} is
indeed $p \vee p'$. This also shows that Lemma \ref{elm} can not
be essentially improved, because  $C_p=R_{p'}\subset S_{p'}$ is a
complemented subspace and hence we can not have cotype $2$, at
most cotype $p$. However, for $p=2$ we know that the exponent is
not attained in general because the little Grothendieck inequality
fails in this form \cite{J2}. Also hyperfiniteness is necessary,
because in the free group algebra $VN(\mathbf{F}_{\infty})$ every
diagonal operator $\Delta_{\mu}(\la(g_k)) = \mu_k \, \la(g_k)$,
$g_k$ the generators, factors completely through
$L_p(VN(\mathbf{F}_{\infty}))$ whenever $\Delta_{\mu}: R_p \cap
C_p \to R \cap C$ is completely bounded. Note here that the span
of the generators is completely complemented (see \cite{P4}) and
we may therefore view these maps as defined on
$VN(\mathbf{F}_{\infty})$. That is, $\mu \in \ell_{2p}$. Hence the
eigenvalues are not in $\ell_p$.}
\end{remark}

\subsection{Applications III. Fourier multipliers}

Our last application is devoted to Fourier multipliers. Let $G$ be
a discrete group and let $VN(G)$ stand for the finite von Neumann
algebra generated by the left regular representation $\la$. Given
a function $\phi: G \to \C$, the corresponding Fourier multiplier
$\lambda(g) \mapsto \phi(g) \lambda(g)$ will be denoted by
$T_\phi$.

\begin{corollary}
If $2 \le p < q < \infty$ and if $$T_{\phi}: VN(G) \to
L_p(VN(G))$$ is completely bounded, then $T_{\phi}: L_q(VN(G)) \to
L_p(VN(G))$ satisfies $$\big\| T_{\phi}: L_q(VN(G)) \to L_p(VN(G))
\big\|_{cb} \le c(p,q) \, \big\| T_{\phi}: VN(G) \to L_p(VN(G))
\big\|_{cb}.$$
\end{corollary}

\dem The algebra $\C[G]$ of finite sums $\sum_g \al_g \la(g)$ is
dense in $L_{p'}(VN(G))$ and $T_{\phi}^*(\C[G]) \subset \C[G]$.
This shows that $T_{\phi}$ is normal. Theorem A gives two norm $1$
elements $a,b \in L_{2q}(VN(G))$ and a cb-map $v: L_q(VN(G)) \to
L_p(VN(G))$ such that $T_{\phi}(x) = v(axb)$. Let $\pi: VN(G) \to
VN(G) \bar{\ten} VN(G)$ be the representation given by
$\pi(\la(g)) = \la(g) \ten \la(g)$. Let us show that the map
$$\Lambda_{ab}: x \in L_q(VN(G)) \mapsto (\mathbf{1} \ten a) \,
\pi(x) \, (\mathbf{1} \ten b) \in L_q(VN(G) \bar{\ten} VN(G))$$ is
completely contractive. This is obvious for $q=\infty$, while for
$q=2$
\begin{eqnarray*}
\Big\| \sum_g \al_g \la(g) \ten a \la(g) b \Big\|_2^2 & = & \summ_g
|\al_g|^2 \big\| a\la(g)b \big\|_2^2
\\ & \le & \|a\|_4^2 \, \|b\|_4^2 \, \summ_g |\al_g|^2 \\ [5pt]
& = & \|a\|_4^2 \, \|b\|_4^2 \, \Big\| \summ_g \al_g \la(g)
\Big\|_2^2.
\end{eqnarray*}
On the other hand, note that $id \ten v: L_q(VN(G) \bar{\ten}
VN(G)) \to L_p(VN(G) \bar{\ten} VN(G))$ is completely bounded.
Indeed, $id \ten v: L_q(L_q) \to L_q(L_p)$ is clearly completely
bounded and the inclusion $L_q(L_p) \subset L_p(L_p)$ is
completely contractive. The latter assertion follows regarding the
involved spaces as conditional $L_p$ spaces and using
interpolation. Combining this with $\Lambda_{ab}$ we find that
$$\pi \big( T_{\phi}(\la(g) \big) = \phi(g) \la(g) \ten \la(g) =
\la(g) \ten v \big( a \la(g) b \big) = (id \ten v) \Lambda_{ab}
(\la(g) \ten \la(g)).$$ Finally, we observe that $\pi: L_p(VN(G))
\to L_p(VN(G)\bar{\ten} VN(G))$ is a completely isometric
embedding. This follows from the $L_p$ version of Fell absorption
principle \cite{PP}. Therefore, we conclude that $T_{\phi} =
\pi^{-1} (id \ten v) \Lambda_{ab}$ is completely bounded. \fin

\noindent \textbf{Acknowledgement.} We would like to thank the
referee for a very careful reading of the paper. His/her
suggestions have led to a much more transparent presentation.

\bibliographystyle{amsplain}

\enlargethispage{2cm}

\vskip10pt

\hfill \noindent \textbf{Marius Junge} \\
\null \hfill Department of Mathematics
\\ \null \hfill University of Illinois at Urbana-Champaign \\
\null \hfill 1409 W. Green St. Urbana, IL 61891. USA \\
\null \hfill\texttt{junge@math.uiuc.edu}

\

\hfill \noindent \textbf{Javier Parcet} \\
\null \hfill Instituto de Ciencias Matem{\'a}ticas \\ \null \hfill
CSIC-UAM-UC3M-UCM \\ \null \hfill Consejo Superior de
Investigaciones Cient{\'\i}ficas \\ \null \hfill Serrano 121.
28006, Madrid. Spain \\ \null \hfill\texttt{javier.parcet@uam.es}
\end{document}